\documentclass{amsart}
\usepackage{amssymb,amscd}

\newtheorem{theorem}{Theorem}[section]
\newtheorem{lemma}[theorem]{Lemma}
\newtheorem{corollary}[theorem]{Corollary}
\newtheorem{proposition}[theorem]{Proposition}
\newtheorem{lemmalg}[theorem]{Lemma/Algorithm}

\newtheorem{propalg}[theorem]{Proposition/Algorithm}

\theoremstyle{definition}
\newtheorem{definition}[theorem]{Definition}
\newtheorem{example}[theorem]{Example}
\newtheorem{remark}[theorem]{Remark} 
\numberwithin{equation}{section}

\newcommand\Complex{\mathbf{C}} 
\newcommand\R{\mathbf{R}}
\newcommand\Q{\mathbf{Q}}
\newcommand\Z{\mathbf{Z}}
\newcommand\q{q} 
\newcommand\F{\mathbf{F}}
\newcommand\Fq{\F_\q}
\newcommand\Fql{\mathbf{F}_{\q^\ell}}

\newcommand\tensor{\otimes}
\newcommand\isomorphic{\cong}
\DeclareMathOperator{\codim}{codim}

\DeclareMathOperator{\Sym}{Sym}
\newcommand\Directsum{\bigoplus}

\newcommand\union{\cup}

\newcommand\intersect{\cap}

\newcommand\abs[1]{{\left|#1\right|}}

\DeclareMathOperator{\image}{image}

\newcommand\Projective{{\bf P}} 
\newcommand\fieldk{k} 
\newcommand\fieldkbar{{\overline{\fieldk}}} 
\newcommand\C{C} 
\newcommand\OC{\mathcal{O}_{\C}} 
\newcommand\SymdC{{\Sym^\degd(\C)}} 
\newcommand\linebundleL{\mathcal{L}} 
\newcommand\linebundleM{\mathcal{M}} 
\newcommand\V{V} 
\newcommand\Vprime{\V'} 
\newcommand\s{s} 
\newcommand\tu{t} 
\newcommand\g{g} 
\newcommand\D{D} 
\newcommand\E{E} 
\newcommand\divF{F} 
\newcommand\degd{d} 
\newcommand\dbar{T} 
\newcommand\HoOC[1]{H^0\bigl(\OC(#1)\bigr)}
\newcommand\HoL{H^0(\linebundleL)} 
\newcommand\HoM{H^0(\linebundleM)} 
\newcommand\pointP{P} 
\newcommand\W{W} 
\newcommand\Wprime{\W'}
\newcommand\WD{\W_\D} 
\newcommand\WprimeD{\Wprime_\D} 
\newcommand\WE{\W_\E} 
\DeclareMathOperator{\flip}{Flip} 
\DeclareMathOperator{\infl}{Infl} 
\DeclareMathOperator{\defl}{Defl} 
\DeclareMathOperator{\addflip}{Addflip} 
\newcommand\HoLm[1]{H^0(\linebundleL( - #1))} 
\newcommand\HoMm[1]{H^0(\linebundleM( - #1))}

\newcommand\mul{\mu} 
\newcommand\HotwoL{H^0(\linebundleL^{\tensor2})} 
\newcommand\HotwoLm[1]{H^0(\linebundleL^{\tensor2}( - #1))} 

\newcommand\igs{IGS}
\newcommand\anigs{an \igs}

\newcommand\eps{\epsilon} 
\newcommand\h{h} 
\newcommand\cwomega{2.376} 
\newcommand\repA{Representation A}
\newcommand\repB{Representation B}
\newcommand\repBzero{Representation B$_0$}

\newcommand\HoB[1]{H^0\bigl(#1\bigr)}

\begin{document}

\title{Asymptotically fast group operations on
Jacobians of general curves\hskip 1.6em}

\author{Kamal Khuri-Makdisi}
\address{Mathematics Department and Center for Advanced Mathematical Sciences,
American University of Beirut, Bliss Street, Beirut, Lebanon}
\email{kmakdisi@aub.edu.lb}
\subjclass[2000]{11Y16, 14Q05, 14H40, 11G20}
\thanks{July 2, 2006}

\begin{abstract}
Let $C$ be a curve of genus $g$ over a field $k$.  We describe
probabilistic algorithms for addition and inversion of the classes of
rational divisors in the Jacobian of $C$.  After a precomputation, which is
done only once for the curve $C$, the algorithms use only linear algebra in
vector spaces of dimension at most $O(g \log g)$, and so take 
$O(g^{3 + \epsilon})$ field operations in $k$, using Gaussian
elimination. Using fast algorithms for the linear algebra, one can improve
this time to $O(g^{2.376})$.  This represents a significant improvement
over the previous record of $O(g^4)$ field operations (also after a
precomputation) for general curves of genus $g$.
\end{abstract}

\maketitle

\section{Introduction}
\label{section1}

Let $\C$ be a smooth projective geometrically irreducible algebraic curve
of genus $\g$ over a field $\fieldk$.  
The Jacobian variety $J$ of $\C$ is a $\g$-dimensional algebraic
group that parametrizes the degree zero divisors on $\C$, up to linear
equivalence.  The Jacobian plays a crucial role both in the theory and in
the applications of the curve $C$, including cryptography and computational
number theory.  For all but the smallest $\g$, it appears impractical to
implement the group $J(\fieldk)$ algorithmically using an
embedding of $J$ into a projective space $\Projective^N$: if we embed $J$
using the complete linear series attached to $3\Theta$ or $4\Theta$
($\Theta$ being the theta divisor), then the equations of $J$ can be
described, but the dimension $N$ grows exponentially with $\g$; on the
other hand, if we use an incomplete linear series, then the equations
defining $J$ become much more complicated.  Instead, algorithms for
$J(\fieldk)$ generally work directly with $\fieldk$-rational divisors on
$\C$, and keep track of their linear equivalence to reduce ``complicated''
divisors to simpler ones as needed.  This gives a computational handle on
the Picard group, $\text{Pic}^0_\fieldk(\C)$, which is a subgroup of
$J(\fieldk)$ (the two groups agree if $C(\fieldk)$ is nonempty).  We shall
nevertheless frequently abuse terminology and refer to the Jacobian
instead of to the Picard group.

In this article, we present what we believe are asymptotically the fastest
algorithms to date that implement the group law on the Picard group of a
\emph{general} curve $\C$, as the genus $\g$ grows.  This assumes that $\C$
is given in one of two specific forms, which we call ``\repA'' and
``\repB,'' with respect to which we can also represent divisors $\D$ on
$\C$.  If we start with equations for $\C$, we need to 
do a single initial precomputation to bring $\C$ into one of these two
forms.  For \repA, this involves computing two Riemann-Roch spaces of the
form $\HoOC{\D_1}$ on $\C$ and a setting up a ``multiplication table''
$\mul$ between them, once and for all.  For \repB, we also need to describe
the values of a basis for a space $\HoOC{\D_1}$ at sufficiently many points
of $\C$; this allows us to speed up the multiplication $\mul$, in a way
analogous to representing polynomials by their values at many points
instead of by their coefficients.  After that, our algorithms boil down to
linear algebra on certain matrices of size $O(\g) \times O(\g \log \g) =
O(\g) \times O(\g^{1+\eps})$, which arise from subspaces of the
Riemann-Roch spaces above.
For \repA, our algorithms attain a complexity of $O(\g^{3+\eps})$ field
operations in $\fieldk$ per group operation (such as addition or negation)
in the Jacobian, and this complexity holds even if we use Gaussian
elimination rather than asymptotically faster algorithms for linear
algebra.  In the case of \repB, the complexity is determined by the linear
algebra.  The current best algorithms~\cite{CoppersmithWinograd} allow us
to attain a complexity of $O(\g^{\cwomega})$ using \repB.
Our algorithms are straightforward to implement and
analyze --- the author had an easy time programming the algorithms for the
Jacobian group in GP/PARI \cite{PARI}, for the case of ``\repA,'' in a
fairly short program file --- but we naturally need more sophisticated
techniques to prove that our algorithms give the correct answer.

Our algorithms are probabilistic, since they have to find certain
intermediate data (``\anigs'' of a divisor $\D$, defined in
Section~\ref{section3r}) for the computation by random search; the above
complexity actually describes the expected number of field operations
needed by our algorithms.  Each trial to find \anigs{} for $\D$ has a
probability of success greater or equal to $1/2$, and 
we can recognize \anigs{} once we have found it, so our algorithms are
guaranteed to terminate with a correct result.  Thus our probabilistic
algorithms are of Las Vegas type.
We have measured complexity by counting field operations in $\fieldk$
instead of, say, bit operations, due to potential ``coefficient explosion''
in $\fieldk$.  This is not an issue if $\fieldk$ is finite, but is
unavoidable if $\fieldk = \Q$ (more generally, for number fields), since
adding points on the Jacobian tends to increase their arithmetic height.
This growth of coefficients will occur even if we carry out our linear
algebra over $\Q$ in the best possible way, for example by incorporating
LLL reduction throughout our algorithms.

Prior to the results of this article, the best algorithms for Jacobians of
general curves had a complexity of $O(\g^4)$ after the initial
precomputations, and were deterministic.  The complexity $O(\g^4)$ was
attained both in the 1999 Ph.D.{} thesis of F. Hess~\cite{HessThesis} (see
also \cite{HessArticle}), and in a 2001 preprint of the author (published
as~\cite{KKM}), whose methods we adapt and extend for this 
article.  The methods of Hess, and of several predecessors of whom we cite
only~\cite{Cantor}, can be called ``arithmetic'': they begin with
a degree $n$ map $\varphi: \C \to \Projective^1$, and view the function
field $\fieldk(\C)$ as a degree $n$ extension of $\fieldk(x)$.  Then
$J(\fieldk)$ is essentially an ideal class group attached to $\fieldk(\C)$,
and we compute with ideals of the integral closure of $\fieldk[x]$ in
$\fieldk(\C)$ by representing them as lattices (i.e., free modules) over
$\fieldk[x]$; one has to also consider the points of $\C$ lying over
$\infty \in \Projective^1$, and the implementation is somewhat involved.
The methods of Hess and his predecessors work best if the minimum
gonality $n = \deg \varphi$ remains bounded as $\g$ grows (for example,
\cite{Cantor} applies only to hyperelliptic curves, for which $n=2$); in
that case, their algorithms generally have complexity $O(\g^2)$.  However,
their methods are sensitive to $n$, and if $n$ grows linearly with $\g$, as
is the case\footnote{By 
\cite{GriffithsHarris}, page 261, a general curve of genus $\g$ over
$\Complex$, or more generally over an algebraically closed field, has
gonality $\left\lfloor (g+1)/2 \right\rfloor + 1$; over $\fieldk$, the
gonality can be higher.  We also
note the result of~\cite{Abramovich} that the gonality of a modular curve
such as $X_0(N)$ also grows linearly with the genus, at least over
$\Complex$; interestingly, our algorithms are particularly suited for
modular curves, since it is easier to describe them using \repA{} or
\repB{} than by finding nice equations.}
for general curves of genus $\g$, then the complexity of Hess' algorithms
rises to $O(\g^4)$, as mentioned above. 

In contrast, the methods in~\cite{KKM} and in the present article can be
called ``geometric,'' in that we work with an embedding $\iota: \C \to
\Projective^n$.  We choose $n$ moderately large, but still $O(\g)$, and the
two Riemann-Roch spaces that we need to compute are the restriction of
linear and quadratic functions from the projective space to $\C$.
The ``multiplication map'' $\mul$ then multiplies two linear functions to
produce a quadratic function.  Once this is in place, the rest is linear
algebra (the reader may wish to compare our approach with another use of
linear algebra to study Jacobians in~\cite{Anderson}, where linear algebra
on Riemann-Roch spaces and invariant theory are used to describe explicit
equations for the Jacobian).
In contrast to our methods, earlier ``geometric'' algorithms for Jacobians
(\cite{HuangIerardi} and \cite{Volcheck}) preferred to work with $n$ as
small as possible, preferably $n=2$, even if this meant using a
singular plane curve birational to $\C$.  Their algorithms involved fairly
elaborate computations with polynomials of degree $O(\g)$, to say nothing
about the problems with singularities. The resulting complexity of those
algorithms was $O(\g^7)$ after precomputations, and so those methods were
superseded by the algorithms of Hess.  The author hopes that this article
and its predecessor~\cite{KKM} will revive interest in the geometric
approach to algorithms for curves.

We also hope that this article will support a point of view explained in
the introduction of~\cite{KKM}, namely, that it is profitable to do
computational algebraic geometry with varieties embedded in Grassmannians.
Here we represent points on Grassmannians as subspaces of a fixed vector
space $\V$, and use linear algebra throughout; we do \emph{not} embed the
Grassmannian variety into projective space, as the ambient projective space
would be too large.  In our setting, we represent a divisor $\D$ of degree
$\degd$ as a codimension $\degd$ subspace $\WD$ of $\V$, which we can
interpret as mapping the symmetric power variety $\SymdC$ into a
Grassmannian.  We take $\degd \geq 2\g$, instead of the more usual approach
$\degd = \g$, because this simplifies our algorithms (essentially since the
fibre in $\SymdC$ over a point of the Jacobian always has the same
structure, a point used notably in Chow's projective
construction of the Jacobian~\cite{Chow}).  We of course include an
algorithm that determines whether two elements of $\SymdC$ represent the
same point on the Jacobian.  For all this and more background, the reader
is encouraged to consult~\cite{KKM} alongside this article.

The speedup in our new algorithms comes partially from the speedup of
multiplication in \repB; however, the most significant improvement is due
to our using \anigs{} for $\D$ instead of the whole space $\WD$ at some
strategic moments.  This allows us to scale down the size 
of the matrices on which we need to do linear algebra, from
$O(\g) \times O(\g^2)$ in~\cite{KKM} to 
$O(\g) \times O(\g^{1+\eps})$ in this article.  It turns out that the
larger matrices of~\cite{KKM} contain redundant data, but it is still not
clear if one can remove the redundant data by a fast deterministic
algorithm.  This is why our algorithms are probabilistic.

The author gratefully thanks the following institutes for their support
during the periods when the main results of this article and its
predecessor were obtained: the Clay Mathematics Institute, for funding
research visits to the U.S.~in the summers of 2000 and 2003; the William
and Flora Hewlett Foundation, for supporting a semester of paid research
leave at the American University of Beirut in fall 2002; and the TEMPUS
program of the European Union, for funding a research visit to France and
Austria in the summer of 2003.  The author also thanks Universit\'e Paris
XIII and Princeton University for their hospitality during the initial
stages of preparing the manuscript.  The author gratefully thanks G. Frey
for an invitation to Essen to lecture on the results of this article in the
summer of 2005, and for many stimulating conversations on that occasion
with him and with C. Diem and F. Hess, to whom the author is also grateful.
Finally, the author thanks the referees of this article, who carefully read
a rather lengthy earlier draft~\cite{KKMextendedpreprint}, and made
extensive comments and suggestions that helped to streamline the
presentation and to produce a much improved and more compact version of
this article.  This version omits a few details and auxiliary results,
particularly in Section~\ref{section5r}, and the reader who wishes more
explanations of those points may wish to read selected portions of the
earlier manuscript.

\begin{remark}
\label{remark1.1}
We have slightly changed notation between this paper and~\cite{KKM}.
We now use ``multiplicative'' notation to refer to line bundles on the
curve $\C$, instead of the ``additive'' notation that was used in most of
the previous paper (actually, the previous paper occasionally used
multiplicative notation as well).  Here is a small table of old vs.~new
notation.
\begin{center}
\begin{tabular}{|l|l|}
\hline
Old Notation: & New Notation: \\
\hline
\hline
$\linebundleL_1 + \linebundleL_2$ & $\linebundleL_1 \tensor \linebundleL_2$ \\
\hline
$H^0(\D_1 - \D_2)$ & $\HoOC{\D_1 - \D_2}$ \\
\hline
$H^0(2\linebundleL - \D - \E)$ & $\HotwoLm{\D - \E}$ \\
\hline
\end{tabular}
\end{center}
\end{remark}

\begin{remark}
\label{remark1.2}
Throughout this article, we will view an $m\times n$ matrix $M$ (always
with entries in $\fieldk$) as a linear transformation from $\fieldk^n$ to
$\fieldk^m$, viewed as column vectors. Thus $v \in \fieldk^n$ is mapped to
$Mv \in \fieldk^m$, and the notations $\ker M$ and $\image M$ should be
interpreted accordingly.  We shall need to refer to the complexity of the
linear algebra steps in our algorithms, which include computing a
kernel or an image of $M$ and/or reduced row and column echelon forms, as
well as multiplying matrices.  (See Chapter~16 of
\cite{BurgisserClausenShokrollahi} or Chapter~6 of \cite{AhoHopcroftUllman}
for a reduction of general linear algebra to matrix multiplication.)
We denote by $\omega$ the smallest exponent such that linear algebra on
square $n\times n$ matrices has complexity $O(n^{\omega+\eps})$, measured
in field operations in $\fieldk$.  The current
record~\cite{CoppersmithWinograd} is $\omega < \cwomega$, and it is 
conjectured that $\omega = 2$.  Gaussian elimination gives $\omega \leq 3$
elementarily, and in fact the complexity of Gaussian elimination on a
rectangular $m\times n$ matrix is $O(mn\min(m,n))$.
\end{remark}

\begin{remark}
\label{remark1.3}
We use in this article the notation
$\lceil a \rceil$ for the ceiling of $a \in \R$, i.e.,
\begin{equation}
\label{equation1.1}
\lceil a \rceil = \min \{ n \in \Z \mid n \geq a \}.
\end{equation}
\end{remark}

\section{Representing the curve, and basic linear algebra operations}
\label{section2r}

In this section, we describe how we represent the curve $\C$ and how we
implement the basic building blocks of our algorithms via linear algebra.
We shall assume that the curve $\C$ comes equipped with a line bundle
$\linebundleL$ of moderately large, but not too large, degree:
\begin{equation}
\label{equation2r.1}
 \deg \linebundleL \geq 2\g+2,
 \qquad\text{ but nonetheless } \quad \deg \linebundleL = O(\g).
\end{equation}
(We will typically take $\deg \linebundleL = 6\g$ in applications.)  
We define $\fieldk$-vector spaces $\V$, $\Vprime$ by
\begin{equation}
\label{equation2r.2}
\V = H^0(\C, \linebundleL) = \HoL, \qquad
   \Vprime = H^0(\C, \linebundleL^{\tensor 2}) = \HotwoL.
\end{equation}
We also introduce the notation for dimensions and degrees
\begin{equation}
\label{equation2r.3}
\begin{split}
\Delta = \deg \linebundleL, &\qquad \qquad
    \Delta' = \deg \linebundleL^{\tensor 2} = 2 \Delta  \\
\delta = \dim \V = \Delta + 1 - \g, &\qquad \qquad
\delta' = \dim \V' = 2\Delta + 1 - \g. \\
\end{split}
\end{equation}
Note that all the above quantities are $O(\g)$.  All our algorithms work
over the field $\fieldk$, but for some of our proofs we need to consider
points of $\C$ and elements of $\V$ defined over the algebraic closure
$\fieldkbar$.

The most important ingredient in our description of $\C$ is then the
multiplication map $\mul$ on global sections,
\begin{equation}
\label{equation2r.4}
\mul: \V \tensor \V = \HoL \tensor \HoL \to \Vprime = \HotwoL.
\end{equation}
We will use the shorthand notation
\begin{equation}
\label{equation2r.5}
\s\cdot \tu = \mul(\s\tensor\tu) \in \Vprime
\qquad \text{ for } \s, \tu \in \V.
\end{equation}
Before specifying the precise form in which we represent $\mul$
algorithmically, we note the following.
\begin{proposition}
\label{proposition2r.1}
We can determine $\C$ and $\linebundleL$ up to isomorphism from a knowledge
of the multiplication map $\mul$.
Moreover, given vector spaces $\V, \Vprime$ and a map $\mul$, it is
possible to determine whether they come from a pair $(\C, \linebundleL)$ as
above.
\end{proposition}
\begin{proof}
For the first statement, assume that we are assured of the existence of
some pair $(\C, \linebundleL)$, but that we only know the map $\mul$.
We claim that the kernel of $\mul$ encodes equations for $\C$.  
Indeed, consider the embedding $\iota:\C \to \Projective(\V) = 
\Projective^{\delta - 1}$ given by $\linebundleL$.  
Since $\Delta \geq 2\g+2$, this embedding is projectively normal (in
particular, $\mul$ is surjective), and the homogeneous ideal
$I_\C \subset \Sym^* \V$ defining $\iota(\C)$ is generated by quadrics (see
for example \cite{Lazarsfeld}).  Concretely, we can identify $\Sym^* \V$
with $\fieldk[T_1, \dots, T_\delta]$, upon choosing a basis $\{T_1, \dots,
T_\delta\}$ of $\V$.  The kernel of $\mul$ trivially contains all
``commutators'' $T_i \tensor T_j - T_j \tensor T_i$.  After we quotient out
by these commutators, the image of $\ker \mul$ inside the symmetric square
$\Sym^2 \V$ then corresponds to the degree~$2$ elements of $I_\C$.
Since these generate $I_\C$, we can hence recover $\C$; we also obtain
$\linebundleL$ as the pullback of $\mathcal{O}_{\Projective(V)}(1)$ to
$\C$.

For the second statement, we check first that $\mul$ is surjective and
symmetric (i.e., $T_i \cdot T_j = T_j \cdot T_i$ for all $i,j$).  The
kernel of $\mul$, when projected to $\Sym^2 \V$, then corresponds to a
space of degree~$2$ polynomials in $\fieldk[T_1, \dots, T_\delta]$, and we
let $I$ be the ideal generated by (a basis for) this space of degree~$2$
polynomials.  We then check that the ideal $I$ is saturated and that it
defines a smooth projective curve $\C$ (e.g., using Gr\"obner bases).  We
then determine the degree $\Delta$ and genus $\g$ of this curve $\C$ from
the Hilbert series of $I$; this again gives $\linebundleL$ as the pullback
of $\mathcal{O}_{\Projective(V)}(1)$, with $\V \subset \HoL$.  We finally
verify that $\dim \V = \Delta + 1 - \g$ to ensure that $\V = \HoL$, i.e.,
that our embedding of $\C$ comes from the \emph{complete} linear series. 
\end{proof}

Note that in light of the above proposition, we can view $\V$ as the space
of linear ``functions'' and $\Vprime$ as the space of
quadratic functions \emph{on the curve $\C$} with respect to the projective
embedding $\iota$.

For algorithmic purposes, we represent our knowledge of $\V$,
$\Vprime$, and $\mul$ in either of two ways, \repA{} and \repB, with the
former more straightforward, and the latter asymptotically faster.  We also
single out a simple special case \repBzero{} of \repB, both for reasons of
exposition and of ease of implementation.  See Example~\ref{example2r.105}
at the end of this section for an example of \repA{} and \repB.

\begin{enumerate}
\item
\textbf{\repA}:
This method works over all fields.  We choose bases $\{T_1, \dots,
T_\delta\}$ for $\V$ and $\{U_1, \dots, U_{\delta'}\}$ for $\Vprime$,
thereby identifying $\V$ and $\Vprime$ with the spaces of column vectors
$\fieldk^\delta$ and $\fieldk^{\delta'}$.  Knowledge of $\mul$ is then
encoded as a multiplication table, i.e., by storing the coefficients
$c_{ijk}$ in each identity
\begin{equation}
\label{equation2r.6}
T_i \cdot T_j = \mul(T_i \tensor T_j) = \sum_k c_{ijk} U_k.
\end{equation}
It is convenient to store this information as a collection 
$\{M_1, \dots, M_\delta\}$ of matrices, each of size $\delta' \times
\delta$, such that $M_i$ describes the linear transformation
``multiplication by $T_i$'' from $\V$ to $\Vprime$:
\begin{equation}
\label{equation2r.7}
M_i = (c_{ijk})_{k,j} =
\begin{pmatrix}
c_{i11} & c_{i21} & \dots & c_{i\delta 1} \\
c_{i12} & c_{i22} & \dots & c_{i\delta 2} \\
\vdots & \vdots & \ddots & \vdots \\
c_{i1\delta'} & c_{i2\delta'} & \dots & c_{i\delta\delta'}\\
\end{pmatrix}.
\end{equation}
\item
\textbf{\repBzero}:
We take a divisor $\D_1$ such that $\linebundleL = \OC(\D_1)$.  We also
assume that we can find $N = \Delta'+1$ distinct points $\pointP_1, \dots,
\pointP_N \in \C(\fieldk)$ that are not in the support of $\D_1$.  This is
a nontrivial assumption if $\fieldk$ is a number field, but is easy to
arrange in cases of interest to cryptography, where $\fieldk$ is a finite
field of large cardinality.  We then represent $\V$ and $\Vprime$ as
certain subspaces of $\fieldk^N$: namely, we have injections of vector
spaces $\V \to \fieldk^N$ and $\Vprime \to \fieldk^N$ given by 
\begin{equation}
\label{equation2r.8}
\s \mapsto (\s(\pointP_1), \dots, \s(\pointP_N)),
\end{equation}
viewing $\s \in \V$ (respectively $\Vprime$) as a meromorphic function on
$\C$ with poles at $\D_1$ (respectively $2\D_1$).  Then the multiplication
map $\mul$ is simply pointwise multiplication, since $\s\cdot \tu$
corresponds to 
$(\s(\pointP_1)\tu(\pointP_1), \dots, \s(\pointP_N)\tu(\pointP_N))$.
(Thus \repBzero{} is analogous to representing a polynomial $f(x) \in
\fieldk[x]$ of bounded degree by its vector of values $(f(a_1), \dots,
f(a_N))$ at sufficiently many points, in order to speed up the
multiplication of polynomials.)  In this setting, we represent $\C$ by our
knowledge of the subspaces of $\fieldk^N$ corresponding to $\V$ and
$\Vprime$.  It is most convenient to store an $N\times \delta$ matrix
$A_\V$ whose columns are a basis of $\V$ (viewed as a subspace of
$\fieldk^N$), as well as the equivalent data of an $(N-\delta)\times N$
matrix $K_\V$ whose kernel is the subspace $\V$.  It turns out not to be
necessary to store a basis for the subspace $\Vprime$, but we can always
recover it, if needed, from the fact that $\mul$ in~\eqref{equation2r.4} is
surjective.  Note that there is no need to store any information that
describes the map $\mul$.
\item
\textbf{\repB}: 
Even if we cannot find enough $\fieldk$-rational points on $\C$, we can
still work with the following generalization under the mild assumption
of~\eqref{equation2r.10} (e.g., it is sufficient to assume that $\fieldk$
is perfect).  We take a $\fieldk$-rational effective divisor $Z$ on $\C$,
of degree $N = \deg Z = O(\g)$, such that $\HotwoLm{Z} = 0$.  (We chose
$N = \Delta' + 1$ and $Z = \pointP_1 + \dots + \pointP_N$ for \repBzero.)
We then wish to represent elements of $\V$ and $\Vprime$ by their
``values'' at the points of $Z$.  Here the values of a global section
$\s \in \V = \HoL$ at $Z$ are given by the image of $\s$ in
$\HoB{\linebundleL_Z}$, where we define the sheaf $\linebundleL_Z =  
\linebundleL/\linebundleL(-Z)$.  We similarly define
$\linebundleL^{\tensor2}_Z =
     \linebundleL^{\tensor2}/\linebundleL^{\tensor2}(-Z)$,
and view the values of an element of  
$\Vprime$ at $Z$ as belonging to $\HoB{\linebundleL^{\tensor2}_Z}$.  By
design, the natural $\fieldk$-linear map $\V \to \HoB{\linebundleL_Z}$ is
injective, and similarly for $\Vprime$.  Moreover, one can find
compatible isomorphisms of sheaves of $\OC$-modules
\begin{equation}
\label{equation2r.9}
\varphi: \linebundleL_Z \isomorphic \mathcal{O}_Z,
\qquad
\varphi^{\tensor2}: \linebundleL^{\tensor2}_Z \isomorphic \mathcal{O}_Z, 
\end{equation}
where $\mathcal{O}_Z = \OC/\OC(-Z)$.  This identifies $\V$ and $\Vprime$ as
$\fieldk$-subspaces of the $N$-dimensional $\fieldk$-algebra $\mathcal{A}
:= \HoB{\mathcal{O}_Z}$, in a way such that the multiplication $\mul$
becomes multiplication in $\mathcal{A}$.  We moreover need to assume the
knowledge of an isomorphism of $\fieldk$-algebras:
\begin{equation}
\label{equation2r.10}
  \mathcal{A} \isomorphic
  \fieldk[x]/(h_1(x)) \times \dots \times \fieldk[x]/(h_r(x)).
\end{equation}
We thus represent elements of $\mathcal{A}$ as tuples of polynomials
$(f_1(x), \dots, f_r(x))$ with $\deg f_i < \deg h_i$.  The coefficients of
the $f_i$ identify $\mathcal{A}$ with $\fieldk^N$ as a $\fieldk$-vector
space; with respect to these coordinates, we can describe $\V$ by matrices
$A_\V$ and $K_\V$ as in the case of \repBzero.  However, in this setting,
we need to carry around the polynomials $h_1(x), \dots, h_r(x)$ in order to
know the multiplication map $\mul$.  Note that multiplying two elements of
$\mathcal{A}$ can be done in time $O(N^{1+\eps}) = O(\g^{1+\eps})$ by
FFT-based methods. 
\end{enumerate}

\begin{remark}
\label{remark2r.3}
It is relatively straightforward to produce \repA{} for a curve that is
given in a more ``classical'' representation.  For instance, we may be
given polynomial equations that describe $\C$ in some projective space
(where the embedding need not be given by a complete linear series).
Alternatively, we may start with a representation of the function field of
$\C$ as an extension $\fieldk(x)[y]$ of the rational function field
$\fieldk(x)$, given by an equation $f(x,y) = 0$; this is tantamount to
choosing a possibly singular plane curve birational to $\C$.  In either of
these two cases, we choose a divisor $\D_1$ of suitably large degree
$\Delta$, and let $\linebundleL = \OC(\D_1)$.  We then use standard
algorithms (\cite{HuangIerardi}, \cite{Volcheck}, \cite{HessArticle}) for
calculating the Riemann-Roch spaces $\V = \HoOC{\D_1}$ and
$\Vprime = \HoOC{2\D_1}$.  The multiplication map $\mul$ is then immediate
in terms of the representation of $\V$ and $\Vprime$ as subsets of the
function field $\fieldk(\C)$.

Another situation where we can produce \repA{} is that of modular curves.
If our curve $\C$ is the completion of a quotient
$\Gamma\backslash\mathcal{H}$ for some congruence subgroup $\Gamma$ acting
on the upper half-plane $\mathcal{H}$, then we do not need to compute
equations for $\C$ directly; instead, we take a suitable weight $n$ (small
values such as $n \in \{2, 3, 4\}$ usually suffice), and let
$\V = \mathcal{M}_n(\Gamma)$ and $\Vprime = \mathcal{M}_{2n}(\Gamma)$ be
the spaces of modular forms of weights $n$ and $2n$ with respect to
$\Gamma$.  The map $\mul$ is then multiplication of modular forms; one way
in which the modular forms can be represented is by their $q$-expansions up
to a suitable order $O(q^N)$, where $N$ is large enough to distinguish
elements of $\Vprime$.  These $q$-expansions can be efficiently computed
using modular symbols (see, e.g., \cite{Stein}).  Note that working with
$q$-expansions is essentially \repB, where the divisor $Z$ is the $N$-fold
multiple of the cusp at infinity.  The author has also investigated \repB{}
for modular curves in the setting where one evaluates the form at several
non-cuspidal points.
\end{remark}

\begin{remark}
\label{remark2r.4}
Given a curve in \repB, one can immediately convert the curve to \repA.
Conversely, given a curve in \repA, we sketch in Section~\ref{section5r}
how to convert this to \repB, under some assumptions on the field
$\fieldk$.
\end{remark}

\begin{remark}
\label{remark2r.2}
For uniformity of notation, we extend the definition of $N$ so that in the
case of \repA, we have $N = \delta$.  Thus both in \repA{} and in \repB, we
will identify $\V$ with a subspace of $\fieldk^N$, viewed as column
vectors:
\begin{enumerate}
\item
If we use \repA, then $\V = \fieldk^N$; in this case we can consider that
$A_\V$ is the $N\times N$ identity matrix.
\item
If we use \repB, then $\V = \image A_\V = \ker K_\V$.
\end{enumerate}
We similarly define $N'$ by $N' = \delta'$ in the case of \repA, and $N' =
N$ in the case of \repB, so that $\Vprime$ is identified with a subspace of
$\fieldk^{N'}$.

We will also need to represent ($\fieldk$-rational) subspaces $\W \subset
\V$ and $\Wprime \subset \Vprime$.  If $r = \dim \W$, then we represent
$\W$ nonuniquely by an $N\times r$ matrix $A_\W$, whose columns give a
basis for $\W$ (viewing the columns as
elements of $\V$).  Thus we have an inclusion $\image A_\W \subset \image
A_\V$ corresponding to the inclusion $\W \subset \V$.  We similarly
represent an $r'$-dimensional subspace $\Wprime \subset \Vprime$ by an $N'
\times r'$ matrix $A_{\Wprime}$ with $\image A_{\Wprime} = \Wprime$.  Note
finally that the numbers $N$ and $N'$, as well as the smaller $r$ and
$r'$, are all $O(\g)$, regardless of whether we use \repA{} or \repB.
\end{remark}

Our algorithms will represent divisors as certain subspaces of $\V$ and of
$\Vprime$, and will all involve the following linear algebra techniques:

\begin{definition}
\label{definition2r.4.5}
Given subspaces $\W \subset \V$ and $\Wprime \subset \Vprime$, and given
elements $\s, \s_1, \dots, \s_\h \in \V$, we define the following:
\begin{enumerate}
\item
The \emph{simple multiplication} $\s \cdot \W$ is the subspace of $\Vprime$
defined by 
\begin{equation}
\label{equation2r.11}
\s \cdot \W = \{\s \cdot \tu \mid \tu \in \W \}.
\end{equation}
\item
The \emph{sum of products} $\s_1 \cdot \W + \dots + \s_\h \cdot \W \subset
\Vprime$ is the usual sum of subspaces.  (We can view this as a ``full
multiplication'' between $S = \text{span}\{\s_1, \dots, \s_\h\}$ and $\W$.)
\item
The \emph{division} $\Wprime \div \{\s_1, \dots, \s_\h\}$ is the
subspace of $\V$ given by
\begin{equation}
\label{equation2r.11.5}
\Wprime \div \{ \s_1, \dots, \s_\h \}
  =  \{ \tu \in \V \mid \tu \cdot \s_i \in \Wprime, 
                    \text{ for all } 1 \leq i \leq \h\}.
\end{equation}
\end{enumerate}
\end{definition}

The above operations were used in the algorithms of~\cite{KKM}, with
$\h = O(\g)$, but we shall only need the case $\h = O(\g^{\eps})$ in this
article.  We can immediately describe the complexity of the above
operations, measured as usual in the number of $\fieldk$-operations.  The
exponent $\omega$ in the complexity of linear algebra was mentioned in
Remark~\ref{remark1.2}.

\begin{propalg}
\label{propalg2r.5}
Assume that $\h = O(\g^{\eps})$.  Using \repA, we can:
\begin{enumerate}
\item
Find one product $\s \cdot \tu$ with complexity $O(\g^3)$.
\item
Compute a simple multiplication $\s \cdot \W$ with complexity $O(\g^3)$.
\item
Compute a sum of products $\s_1 \cdot \W + \dots + \s_\h \cdot \W$ with
complexity $O(\g^{3+\eps})$.
\item
Compute a division $\Wprime \div \{\s_1, \dots, \s_\h\}$ with complexity
$O(\g^{3+\eps})$.
\end{enumerate}
\end{propalg}
\begin{proof}
\begin{enumerate}
\item
Our representation of elements of $\V$ as tuples in $\fieldk^N$, via the
basis $\{T_i\}$ for $\V$, means that we are given $\s = \sum_i
c_i T_i$ in the form of the column vector ${}^t(c_1,\dots,c_N)$.
It is useful to produce the $N'\times N$ matrix $M_\s$ which describes
the linear transformation ``multiplication by $\s$'' from $\V$ to
$\Vprime$:
\begin{equation}
\label{equation2r.12}
M_\s = \sum_i c_i M_i, \qquad\qquad
   M_i \text{ as in \eqref{equation2r.7}}.
\end{equation}
Also viewing $\tu$ as a column vector in $\fieldk^N$, we then compute $\s
\cdot \tu = M_\s \tu$.  Here computing $M_\s$ has complexity $O(\g^3)$,
and multiplying $M_\s \tu$ has complexity $O(\g^2)$.  (Alternatively, we
could have expanded $\s\cdot \tu$ using the coefficients $c_{ijk}$ of the
multiplication table~\eqref{equation2r.6}, for the same complexity.)
\item
We are given the matrix $A_\W$, as in Remark~\ref{remark2r.2}.  Compute the
matrix $M_s$ as above, with complexity $O(\g^3)$; then form the matrix
product $A_{\s \cdot \W} = M_s A_\W$.  We remain within complexity
$O(\g^3)$, even if we use fast matrix multiplication.  Note that the naive
method of multiplying $\s$ by each column of $\W$ would have had complexity
$O(\g^4)$.
\item
Compute the matrices $A_{\s_1 \cdot \W}, \dots, A_{\s_\h \cdot \W}$.  So
far, this requires a complexity of $O(\g^3\h)$.  Then our desired result is
the image of the block matrix 
$A' = \Bigl(A_{\s_1 \cdot \W} \quad \dots \quad A_{\s_\h \cdot \W}\Bigr)$,
whose size is $O(\g) \times O(\g\h)$.  We then find a basis for
$\image A'$ by linear algebra, with complexity $O(\g^3\h)$ if we use
Gaussian elimination, and $O((\g\h)^{\omega+\eps})$ by fast methods.  Our
total complexity is then $O(\g^{3+\eps})$.
\item
Let $r' = \dim \Wprime$.  From the $N'\times r'$ matrix $A_{\Wprime}$, we
use linear algebra to produce an $(N'-r')\times N'$ matrix $K_{\Wprime}$
whose kernel is $\Wprime$; the complexity of this is dominated by what
comes next.  Then our desired result is 
\begin{equation}
\label{equation2r.12.5}
A_{(\Wprime \div \{\s_1, \dots, \s_\h\})} = \ker P, \qquad \text{where }
P = \begin{pmatrix}
K_{\Wprime} M_{\s_1}\\
\vdots \\
K_{\Wprime} M_{\s_\h} \\
\end{pmatrix}.
\end{equation}
This takes complexity $O(\g^3\h)$ to produce the $\{M_{\s_i}\}$, then
$O(\g^{\omega+\eps}\h)$ to obtain $P$.  The matrix $P$ has size
$\bigl((N'-r')\h\bigr) \times N = O(\g\h) \times O(\g)$, and finding its
kernel has a complexity of $O((\g\h)^{\omega+\eps})$ (even if we use
Gaussian elimination, the time is still dominated by finding the
$M_{\s_i}$). 
\end{enumerate}
\end{proof}

Note that for \repA, there is no asymptotic advantage to using fast linear
algebra; we can carry out the operations of
Proposition/Algorithm~\ref{propalg2r.5} using Gaussian elimination with the
same complexity, albeit with a higher implied constant in the  $O(\cdot)$
notation.  On the other hand, \repB{} benefits significantly from fast
linear algebra.

\begin{propalg}
\label{propalg2r.5.3}
Assume that $\h = O(\g^{\eps})$.  Using \repB, we can:
\begin{enumerate}
\item
Find one product $\s \cdot \tu$ with complexity $O(\g^{1+\eps})$.
\item
Compute a simple multiplication $\s \cdot \W$ with complexity $O(\g^{2+\eps})$.
\item
Compute a sum of products $\s_1 \cdot \W + \dots + \s_\h \cdot \W$ with
complexity $O(\g^{\omega+\eps})$.
\item
Compute a division $\Wprime \div \{\s_1, \dots, \s_\h\}$ with complexity
$O(\g^{\omega+\eps})$.
\end{enumerate}
\end{propalg}
\begin{proof}
This is largely the same as the previous result, except that the bottleneck
caused by finding matrices of the form $M_{\s}$ can be bypassed.  We
indicate the necessary modifications.  Note that if we use \repBzero, then
the first two statements hold without including $\eps$ in the exponents.
\begin{enumerate}
\item
Recall that we represent $\s, \tu$ as elements of the algebra
$\mathcal{A}$ (which is just $\fieldk\times \dots \times \fieldk$ for
\repBzero, in which case the result is even easier), and we can multiply
two elements of $\mathcal{A}$ by FFT-techniques.
\item
Either multiply $\s$ by each column of $A_\W$ separately.  We note for
later use the fact that the $N'\times N = N \times N$ matrix $M_\s$ is
block diagonal with a structure that allows fast multiplication by FFT ---
the matrix $M_\s$ is furthermore genuinely diagonal in the case of
\repBzero. Hence the multiplication $M_\s A_\W$ can be done with complexity
$O(\g^{2+\eps})$.  If we want, we can actually produce $M_\s$ by
directly multiplying $\s$ by each element in our basis for $\mathcal{A}
\isomorphic \fieldk^N$.  This also has complexity $O(\g^{2+\eps})$;
it corresponds to replacing $A_\W$ by the identity matrix.
\item
Here it only takes us complexity $O(\g^{2+\eps}\h)$ to produce the matrix
$A'$, so the result follows.
\item
First note that the matrix $P$ must be replaced by a slightly larger
matrix $Q$ that includes an extra subblock $K_\V$ as mentioned in our
descriptions of \repBzero{} and \repB:
\begin{equation}
\label{equation2r.13}
A_{(\Wprime \div \{\s_1, \dots, \s_\h\})} = \ker Q, \qquad \text{where }
Q = \begin{pmatrix} K_\V \\ 
  K_{\Wprime} M_{\s_1}\\
  \vdots \\
  K_{\Wprime} M_{\s_\h} \\
  \end{pmatrix}.
\end{equation}
This ensures that elements of $\ker Q$ genuinely belong to $\V$, which is a
proper subspace of $\fieldk^N$.  This does not affect the asymptotics of
the linear algebra to find $\ker Q$, since $Q$ still has size $O(\g\h)
\times O(\g)$.  As for finding $Q$ in the first place, note that the
product matrices $\{K_{\Wprime} M_{\s_i}\}$ can be computed with complexity
$O(\g^{2+\eps})$.  This is particularly clear for \repBzero, since
$M_{\s_i}$ is a diagonal matrix.  The proof in general uses the
transposition principle.
Indeed, since the complexity using FFT-based algorithms of
multiplying $M_{\s_i} v$ for any column vector $v \in \fieldk^N$ is
$O(\g^{1+\eps})$, it follows that one can just as quickly (perhaps with
a ``larger'' $\eps$) multiply $w M_{\s_i}$ for any $N$-dimensional
\emph{row vector} $w$.  Applying this to the rows of $K_{\Wprime}$, we
obtain our result.  Alternatively, we can give a more pedestrian approach
to finding $Q$; this takes a slightly higher complexity of
$O(\g^{\omega+\eps})$, but does not affect the final complexity of division.
Simply produce all the matrices $M_{\s_i}$, which requires complexity
$O(\g^{2+\eps}\h)$, and then multiply them by a fast algorithm with the
matrix $K_{\Wprime}$.
\end{enumerate}
\end{proof}

All our later algorithms will be built up from the operations that we have
introduced in the above two Proposition/Algorithms \ref{propalg2r.5}
and~\ref{propalg2r.5.3}.  We shall use the following terminology.

\begin{definition}
\label{definition2r.5.6}
A \emph{fast algorithm} is one that requires a complexity of
$O(\g^{3+\eps})$ field operations in $\fieldk$ using \repA, and that
requires a complexity of $O(\g^{\omega+\eps})$ using \repB.  We will also
define \emph{fast probabilistic algorithms} of Las Vegas type to be those
whose expected running time is of the above complexity.  (Recall that a
probabilistic algorithm is called of Las Vegas type if it either returns an
answer which is guaranteed to be correct, with a probability that is
bounded below by a fixed positive number, or if it returns ``failure.''
This is in contrast to Monte Carlo probabilistic algorithms, for which the
answer in the first instance may be wrong, also with a bound on the
probability of error.)
\end{definition}

We conclude this section with a concrete example of a curve as given in
\repA{} and \repB, in order to clarify the precise input to our algorithms.

\begin{example}
\label{example2r.105}
Let $\C$ be the elliptic curve given by the Weierstrass equation
$y^2 = x^3 + 1$ over a field $\fieldk$ not of characteristic $2$ or~$3$.
We choose as our line bundle $\linebundleL = \OC(4P_\infty)$, where
$P_\infty \in \C(\fieldk)$ is the point at infinity.  We choose bases for
$\V$ and $\Vprime$ (which we view as subsets of $\fieldk(\C)$):
\begin{equation}
\label{equation2r.101}
\begin{split}
\{T_1, \dots, T_4\} &= \{1, x, y, x^2\}, \\
\{U_1, \dots, U_8\} &= \{1, x, y, x^2, xy, x^3, x^2y, x^4\}. \\
\end{split}
\end{equation}
Thus, using \repA, we would have $T_2 \cdot T_3 = U_5$ and
$T_3 \cdot T_3 = U_1 + U_6$.  The reader is encouraged to write down the
matrices $M_i$ of~\eqref{equation2r.7}, which will be the entire
description of our curve $\C$; in particular, our representation never
works with the variables $x$ and $y$, but only with the multiplication
table giving each $T_i \cdot T_j$ in terms of the $U_k$'s.

To illustrate \repB, we take $\fieldk=\Q$, and take the divisor $Z$ of
degree $N=9$ to be
\begin{equation}
\label{equation2r.102}
Z = (0,1) + (-1,0) + (2,3) + (2,-3)
    + (2+\sqrt{2},5+4\sqrt{2}) + (2-\sqrt{2},5-4\sqrt{2})
    + 3P_\infty.
\end{equation}
Note that the individual points need not be defined over $\Q$, but the
divisor $Z$ is nonetheless rational over $\Q$. Here we have chosen the map
$\varphi$ of~\eqref{equation2r.9} to be multiplication by $x^{-2}$ at
$P_\infty$ and to be the identity away from $P_\infty$.  In other words,
the natural 
trivialization of $\linebundleL = \OC(4P_\infty)$ on the complement of
$P_\infty$ allows us to directly evaluate elements of $\V$ or $\Vprime$,
viewed as elements of the function field, at the six ``finite'' points of
$Z$; since the values of a $\Q$-rational element at the points
$(2 \pm \sqrt{2}, 5 \pm 4 \sqrt{2})$ are conjugate elements of the extension
$\Q[\sqrt{2}]$, the values at these two points are completely described by a
single element of $\Q[\sqrt{2}]$.  This is equivalent to noting that these
two conjugate points on $\C(\fieldkbar)$ correspond to a single point on
the scheme $\C$, with residue field $\Q[\sqrt{2}]$.  

As for evaluating at the remaining point
$P_\infty$ (to third order), we ``evaluate'' an element $s \in \V$ by
evaluating the function field element $sx^{-2}$, which is regular at 
$P_\infty$, to third order at that point.  More precisely, we take the
first three terms $sx^{-2} = a_0 + a_1 t + a_2 t^2 + O(t^3)$ in the power
series expansion of $sx^{-2}$ in terms of a uniformizer $t$ of the discrete
valuation at $P_\infty$.  (Specifically, we choose $t = x/y$, so that $x =
t^{-2} + O(t^4)$ and $y = t^{-3} + O(t^3)$.  Also, if we wanted to evaluate
an element $s'\in \Vprime$ at $3P_\infty$, we would need to take the
third-order expansion of $s'x^{-4}$ in terms of $t$.)  Putting all this
together, we see that the algebra $\mathcal{A}$ of ``values at $Z$'' can be
identified with
\begin{equation}
\label{equation2r.103}
\mathcal{A} \isomorphic
   \Q \times \Q \times \Q \times \Q \times \Q[u]/(u^2 - 2) \times \Q[t]/(t^3),
\end{equation}
where $u$ corresponds to $\sqrt{2}$, and the ``values'' of the basis
elements of $\V$ at $Z$ are
\begin{equation}
\label{equation2r.104}
\begin{split}
&(1,1,1,1,1+0u,0+0t+0t^2) \in \mathcal{A} \qquad 
      \text{corresponding to } T_1 \leftrightarrow 1,\\
&(0,-1,2,2,2+u,0+0t+t^2) \qquad\qquad
      \text{corresponding to } T_2 \leftrightarrow x,\\
&(1,0,3,-3,5+4u,0+t+0t^2) \qquad\qquad
      \text{corresponding to } T_3 \leftrightarrow y,\\
&(0,1,4,4,12+8u,1+0t+0t^2) \qquad\qquad
      \text{corresponding to } T_4 \leftrightarrow x^2.\\
\end{split}
\end{equation}
Each element of $\mathcal{A}$ above corresponds to a column of the
$9\times4$ matrix $A_\V$; for example, the third column is
${}^t(1,0,3,-3,5,4,0,1,0)$.  The matrix $A_\V$, along with the
identification of $\mathcal{A}$ with $\Q^9$ via~\eqref{equation2r.103}
(especially the polynomial equations $u^2-2=0$ and $t^3 = 0$), then
constitute our description of $\C$ in \repB.  Note that 
we have not bothered to slavishly follow~\eqref{equation2r.10} in the sense
of writing the first four factors of $\mathcal{A}$ as quotients of univariate
polynomial rings instead of as $\Q$ (e.g., by having the first four factors
be $\Q[w]/(w)$ instead).  What we have done instead is to combine ideas
from \repBzero{} and \repB.
\end{example}

\section{Representing divisors; algorithms for divisor classes}
\label{section3r}

We now turn to the representation of divisors on $\C$.  We begin with
some notation.  Given a divisor $\D$ and a $\pointP \in \C(\fieldkbar)$, we
write $v_\pointP(\D)$ for the multiplicity of $\pointP$ in $\D$; hence
$\D = \sum_\pointP v_\pointP(\D) \pointP$, a finite sum.  We write
$(s)_\linebundleL$, or $(s)$ if $\linebundleL$ is understood, for the
divisor of zeros of a nonzero section $s \in \HoL$:
\begin{equation}
\label{equation2r.15}
(\s) = (\s)_\linebundleL
     = \sum_{\pointP \in \C(\fieldkbar)} v_{\linebundleL,\pointP}(\s) \pointP.
\end{equation}
Here $v_{\linebundleL,\pointP}(\s)$ is the valuation of $\s$ at the point
$\pointP \in \C(\fieldkbar)$.  Note that $(\s)$ is an effective divisor,
with $\deg (\s) = \deg \linebundleL = \Delta$.  Moreover, the linear
equivalence class of $(\s)$ is the same as that of the line bundle
$\linebundleL$, and so is independent of the choice of $\s$.  Note also
that since $\s \in \V$ is rational over $\fieldk$, so is the divisor
$(\s)$, even though the individual points where $(\s)$ vanishes might be
defined over an extension of $\fieldk$.

\begin{definition}
\label{definition2r.6}
Let $\D$ be a $\fieldk$-rational effective divisor on $\C$.
\begin{enumerate}
\item
We define the ($\fieldk$-rational) subspaces
\begin{equation}
\label{equation2r.16}
\begin{split}
\WD &= \{ \s \in \V \mid \forall \pointP \in \C(\fieldkbar),
  v_{\linebundleL,\pointP}(\s) \geq v_{\pointP}(\D) \} = \HoLm{\D} \subset
  \V,\\
\WprimeD &= \{ \s' \in \Vprime \mid \forall \pointP \in \C(\fieldkbar),
  v_{\linebundleL^{\tensor2},\pointP}(\s') \geq v_{\pointP}(\D) \} =
\HotwoLm{\D} \subset \Vprime.\\
\end{split}
\end{equation}
Thus $\WD$ and $\WprimeD$ consist respectively of those linear or quadratic
functions on $\C$ that vanish at~$\D$, counting multiplicity.  We allow $\D
= 0$, in which case $\WD = \V, \WprimeD = \Vprime$.
\item
Take a subset $S \subset \V$ containing at least one nonzero element.
We say that $S$ is an \emph{ideal generating set} (abbreviated to \igs) for
$\D$, or equivalently that $\D$ is the \emph{divisor of common zeros} of
$S$, if
\begin{equation}
\label{equation2r.17}
\forall \pointP \in \C(\fieldkbar), \qquad v_{\pointP}(\D)
    = \min\{ v_{\pointP}(\s) \mid \s \in S\}.
\end{equation}
We occasionally abuse terminology and call $S$ \anigs{} for $\WD$.
\end{enumerate}
\end{definition}

Note that the divisor of common zeros of $S$ is the same as that of the
$\fieldk$-subspace of $\V$ spanned by $S$.  The terminology \igs{} comes
from the interpretation of a divisor $\D$ on (an affine part of) $\C$ as an
ideal in a Dedekind domain.

Clearly, \anigs{} for $\D$ exists if and only if the line bundle
$\linebundleL(-\D)$ is base point free, in which case $\WD$ itself (or even
just a basis for $\WD$) will be \anigs.  The divisor $\D$ is then uniquely
determined by any \igs{} $S$, as it can be viewed as the GCD of the
divisors $\{(\s) \mid 0 \neq \s \in S\}$.  Thus we represent our divisors
as follows:

\begin{definition}
\label{definition2r.7}
Assume that $\D$ is an effective $\fieldk$-rational divisor.
By abuse of terminology, we say that $\WD$ is \emph{base point free} if
the line bundle $\linebundleL(-\D)$ is base point free.
\begin{enumerate}
\item
If $\WD$ is base point free, then a \emph{full representation} of $\D$ is
any matrix $A_{\WD}$ whose columns (as in Remark~\ref{remark2r.2}) are a
basis for the subspace $\WD$.
\item
If $\WD$ is base point free, then a \emph{brief representation} of $\D$ is
any \igs{} $\{\s_1, \dots, \s_\h\}$ for $\D$, where we store the
$\s_i \in \V$ as column vectors in $\fieldk^N$.
\end{enumerate}
\end{definition}

In particular, if $\W$ is a subspace of $\V$ whose divisor of common zeros
is $\D$, then any basis for $\W$ can be viewed as a brief representation of
$\D$.  The following proposition collects some elementary facts that play
an important role in our algorithms.

\begin{proposition}
\label{proposition2r.8}
Let $\D$ be an effective $\fieldk$-rational divisor of degree $\degd$
(we allow $\D = 0$).  Recall that
$\Delta = \deg \linebundleL \geq 2\g + 2$. 
\begin{enumerate}
\item
If $\degd \leq \Delta - 2\g$, then $\WD$ is base point free.  Further,
$\dim \WD = \delta - \degd$ has codimension $\degd$ in $\V$.
\item
If $\degd \leq 2\Delta - 2\g$, then a similar statement holds for the
subspace $\WprimeD \subset \Vprime$.
\item
Take a nonzero $\s \in \V$ with $(\s)_\linebundleL = \E$.  Then the simple
multiplication $\s \cdot \WD$ is
\begin{equation}
\label{equation2r.18}
\s \cdot \WD = \Wprime_{\D + \E}.
\end{equation}
If furthermore
$\degd \leq \Delta -2\g$, then both $\WD$ and $\Wprime_{\D+\E}$ are base
point free.
\item
Let $S = \{\s_1, \dots, \s_\h\}$ be \anigs{} for $\D$.  Let $\E$ be an
effective $\fieldk$-rational divisor, preferably but not necessarily such
that $\Wprime_{\D+\E}$ is base point free.  Then the division
$\Wprime_{\D+\E} \div S$ is
\begin{equation}
\label{equation2r.19}
\Wprime_{\D+\E} \div S = \WE.
\end{equation}
\end{enumerate}
\end{proposition}
\begin{proof}
Easy considerations about valuations and the Riemann-Roch theorem;
the main ideas are present in~\cite{KKM}.  Incidentally, one can also
define $\Wprime_\divF \div S$ for arbitrary divisors $\divF$; the result is
then $\W_{\divF\backslash\D}$, in the sense of Proposition/Algorithm 3.9
of~\cite{KKM}.
\end{proof}

Our next goal is to explain that, with good probability, a random selection
of relatively few elements of a base point free space $\WD$ is \anigs{} for
$\D$.  Moreover, it is easy to test whether any given subset of $\WD$
is \anigs, in the setting of our application.  This enables us to convert
easily between the full and brief representations of $\D$.  We first
clarify what we mean by a random selection of elements of $\WD$, and then
state our result precisely.

\begin{definition}
\label{definition2r.9}
Let $\Sigma \subset \fieldk$ be a finite subset, and let $\abs{\Sigma}$ be
its cardinality.  (If $\fieldk$ is itself finite, we usually take $\Sigma =
\fieldk$.)  Let $\W \subset \V$ be a subspace, and
choose once and for all a basis $\{w_1, \dots, w_r\}$ for $\W$.  We define
a \emph{$\Sigma$-random element $\tu \in \W$} to be an element of the form
\begin{equation}
\label{equation2r.20}
\tu = c_1 w_1 + \dots + c_r w_r, \qquad
c_1, \dots, c_r \in \Sigma,
\end{equation}
where the $c_i$ are chosen independently and randomly with respect to the
uniform probability distribution on $\Sigma$.  Our notation does not
indicate the dependence  on the choice of basis $\{w_1, \dots, w_r\}$,
even though this affects the distribution, because the final results on
random selection of \anigs{} are independent of this choice of basis.  Note
that choosing a $\Sigma$-random element $\tu$ requires
$O(r\log \abs{\Sigma})$ random bits to produce $c_1, \dots, c_r$, followed
by $O(r\g) = O(\g^2)$ field operations in $\fieldk$ for the linear
combination.  We will mainly consider sets $\Sigma$ that are not too large:
$\abs{\Sigma} = O(\g)$; it is also reasonable to take $\abs{\Sigma} =
O(1)$, which is the case if $\fieldk$ is a finite field.
\end{definition}

\begin{theorem}
\label{theorem2r.10}
Let $\D$ be an effective $\fieldk$-rational divisor with
$\degd = \deg \D \leq \Delta - 2\g$.
\begin{enumerate}
\item
Take a finite set $\Sigma \subset \fieldk$ as above.  Define
\begin{equation}
\label{equation2r.21}
\h = 1 + \left\lceil 
            \log 2(\Delta - \degd) / \log \abs{\Sigma}
      \right\rceil.
\end{equation}
Take any nonzero $\s_1 \in \WD$, and choose, $\Sigma$-randomly and
independently, $\h - 1$ elements $\s_2, \dots, \s_\h \in \WD$.  Then with
probability greater than or equal to $1/2$, the set
$\{\s_1, \dots, \s_\h\}$ is \anigs{} for $\D$. 
\item
Independently of part~1, assume that $2\g - 1 \leq \degd \leq \Delta$.  Let
$\h$ be any integer, and take elements $\s_1, \dots, \s_\h \in \WD$.  Then 
$\{ \s_1, \dots, \s_\h \}$ is \anigs{} for $\D$ if and only if the sum of
products $\s_1 \cdot \V + \dots + \s_\h \cdot \V$ satisfies
\begin{equation}
\label{equation2r.22}
\s_1 \cdot \V + \dots + \s_\h \cdot \V = \WprimeD.
\end{equation}
\end{enumerate}
\end{theorem}
\begin{proof}
Part~1 follows from Proposition~\ref{proposition4r.3}, with $\linebundleM =
\linebundleL(-\D)$ and $\eta = 1/2$.  Note that the result still holds even
if we choose $\s_2, \dots, \s_\h$ independently and $\Sigma$-randomly from
a subspace $\W \subset \WD$ whose divisor of common zeros is $\D$.
Part~2 is Proposition~\ref{proposition4r.13}.
\end{proof}

\begin{remark}
\label{remark2r.11}
Since both $\Delta$ and $\degd$ are of size $O(\g)$, we therefore can
obtain a randomly chosen \igs{} of size
$\h = O\big(1 + (\log g / \log \abs{\Sigma})\bigr) = O(\g^{\eps})$ in fewer
than two attempts on average.  This is a considerable improvement over
using a basis of $\WD$, which would contain $O(g)$ elements, and which
would slow down the algorithms of Proposition/Algorithms \ref{propalg2r.5}
and~\ref{propalg2r.5.3}.  This (along with the insight to use
\repB) is the source of the essential speedup in this article, compared to
the algorithms of~\cite{KKM}.
\end{remark}

Using the framework of Section~\ref{section2r} and this section, we now
describe how to convert between the full and brief representations of a
divisor $\D$.  We also introduce the important ``flipping'' algorithm.

\begin{propalg}[Deflation]
\label{propalg3r.7}
Assume given a subspace $\WD \subset \V$ which is the full representation
of a divisor $\D$ with $2\g-1 \leq \deg \D \leq \Delta - 2\g$.  Then there
exists a fast probabilistic Las Vegas algorithm that computes a brief
representation $\{\s_1, \dots, \s_\h\}$ of $\D$, with $\h = O(\g^\eps)$.
We call this a \emph{deflation} of $\D$; even though the deflation is not
unique, we still write
\begin{equation}
\label{equation3r.9}
\defl(\WD) = \{ \s_1, \dots, \s_\h \}, \qquad
   \text{where } \{\s_1, \dots, \s_\h\} \text{ is any \igs{} for } \D.
\end{equation}
\end{propalg}
\begin{proof}
We know that $\deg \D = \dim \V - \dim \WD$.  This means that we know the
dimension $\dim \WprimeD = \dim \Vprime - \deg \D$, even though we have not
yet computed the subspace $\WprimeD$.  We now run the following algorithm:
\begin{enumerate}
\item
Compute the value of $\h$ from~\eqref{equation2r.21}, and randomly choose
$\s_1, \dots, \s_\h \in \WD$ as in Theorem~\ref{theorem2r.10}
above.
\item
Form the sum of products $\Wprime = \s_1 \cdot \V + \dots + \s_\h \cdot
\V$ by our fast algorithm.  If $\dim \Wprime \neq \dim
\WprimeD$, then our choice of $\{\s_1, \dots, \s_\h\}$ was not \anigs, so
return to step~1.  Once the $\dim \Wprime = \dim \WprimeD$, stop and output
the $\{\s_i\}$.
\end{enumerate}
The complexity of step~1 (including generating the random bits and
forming each $\s_i$) is $O(\g^2 \h) = O(\g^{2+\eps})$, which can be brought
down slightly if one views producing $\{\s_2, \dots, \s_\h\}$ as a matrix
multiplication of $A_{\WD}$ by a random matrix with entries in $\Sigma$.
As for step~2, we have $\Wprime \subset \WprimeD$, so checking the
criterion of~\eqref{equation2r.22} amounts to comparing dimensions.  Our
choice of the $\{\s_i\}$ fails this test with probability at least~$1/2$,
so the expected number of times that we go through the loop is at
most~$2$.
\end{proof}

Converting back from a brief to a full representation of a divisor, which
we call ``inflation,'' requires \anigs{} for $\V$.  This should be computed
once and for all as part of our precomputations when we store the
representation of $\C$ and $\mul$ for our algorithms.  The rest of our
algorithms do not use inflation, but we include it for completeness.
As for the \igs{} for $\V$, we do not need it to implement the group
operations on divisor classes on $\C$, but we do need to have it available
for the ``membership test'' of Section~\ref{section4r},
which tests whether a given subspace $\W \subset \V$ is equal to some
$\WD$.

\begin{lemmalg}[\igs{} for $\V$]
\label{lemmalg3r.8}
There exists a polynomial-complexity, but not ``fast,'' Las Vegas algorithm
that can be done exactly once as a precomputation to produce \anigs{} for
$\V$.  We shall call the (nonunique) result $\defl(\V)$.
\end{lemmalg}
\begin{proof}
As we wish to produce \anigs{} for the empty divisor $\D = 0$, we cannot
use part~2 of Theorem~\ref{theorem2r.10} here.  We need to go beyond the
linear and quadratic spaces $\V$ and $\Vprime$ to a ``cubic'' space 
$\V'' = H^0(\linebundleL^{\tensor3})$.  Write the product of $\s \in \V$
and $\tu' \in \Vprime$ as $\s * \tu' \in \V''$; then the condition for
$\{\s_1, \dots, \s_\h\} \subset \V$ to be \anigs{} for $\V$ is
\begin{equation}
\label{equation3r.10}
\s_1 * \Vprime + \dots + \s_\h * \Vprime = \V''.
\end{equation}
There is no problem in choosing the $\{\s_i\}$ from $\V$ that have a
probability of at least~$1/2$ of being \anigs{} for $\V$.  Carrying
out the modified sum of products in~\eqref{equation3r.10}, however, needs a
knowledge of the space $\V''$ and of the higher multiplication map
$*: \V \times \Vprime \to \V''$; the problem is to produce this data, after
which checking~\eqref{equation3r.10} is easy.  (The data giving $\V''$ and
$*$ can incidentally be discarded once we find \anigs{} for $\V$.)  To find
this data, we can use \repA{} by Remark~\ref{remark2r.4}.  Then, as in
Proposition~\ref{proposition2r.1}, we let $\{T_1, \dots T_\delta\}$ be a
basis for $\V$, and work with the polynomial algebra
$\fieldk[T_1, \dots, T_\delta]$.  The kernel of $\mul$ allows us to find
generators of the ideal $I_\C$, and we can identify $\V$, $\Vprime$, and
$\V''$ respectively as the portions of the graded algebra
$\fieldk[T_1, \dots, T_\delta]/I_\C$ in degrees $1$, $2$, and~$3$, with the
obvious multiplications.  Thus finding $\V''$ and $*$ can be done by
Gr\"obner bases; the computations involve only linear algebra in the spaces
of polynomials in $\fieldk[T_1, \dots, T_\delta]$ of degree at most~$3$,
whose dimension is $O(\g^3)$.  Thus the computation can be done with a
complexity that is polynomial in $\g$.
\end{proof}

\begin{propalg}[Inflation]
\label{propalg3r.9}
Given a precomputed \igs{} for $\V$, assume we are given a brief
representation $\{\s_1, \dots, \s_\h\}$ of a divisor $\D$, with $\h =
O(\g^\eps)$.  Assume that we know that $\deg \D \geq 2\g - 1$.  Then there
exists a (deterministic) fast algorithm to find the full representation
$\WD$, which we call the \emph{inflation} of the \igs{}
$\{\s_1, \dots, \s_\h\}$:
\begin{equation}
\label{equation3r.11}
\infl(\{\s_1, \dots, \s_\h\}) = \WD, \qquad
  \D = \text{ the divisor of common zeros of }
      \{\s_1, \dots, \s_\h\}.
\end{equation}
\end{propalg}
\begin{proof}
The obvious algorithm is:
\begin{enumerate}
\item 
Calculate the sum of products $\WprimeD = \s_1 \cdot \V + \dots + \s_\h
\cdot \V$.
\item
Use the previously computed \igs{}, $\defl(\V)$, to find $\WD = \WprimeD
\div \defl(\V)$.
\end{enumerate}
\end{proof}

The next Proposition/Algorithm is fundamental for our algorithms on
divisors and divisor classes.  Given $\D$, it allows us to find a
complementary (effective) divisor $\tilde{D}$ such that $\D+\tilde{\D}$ is
in the linear 
equivalence class of $\linebundleL$.

\begin{propalg}[Flipping]
\label{propalg3r.10}
Assume given $\WD$, where $2\g - 1 \leq \deg \D \leq \Delta - 2\g$.
Take a nonzero $\s \in \WD$, and write the divisor of $\s$ as
$(\s)_\linebundleL = \D + \tilde{\D}$.  Then there exists a fast Las Vegas
algorithm to compute the \emph{flip}, $\W_{\tilde{\D}}$, of our divisor:
\begin{equation}
\label{equation3r.12}
\flip(\WD,\s) = \W_{\tilde{\D}}.
\end{equation}
\end{propalg}
\begin{proof}
Compute $\W_{\tilde{\D}} = (\s \cdot \V) \div \defl(\WD)$.  This
works because $\s\cdot \V = \Wprime_{\D + \tilde{\D}}$.
\end{proof}

\begin{remark}
\label{remark3r.11}
We will write $\W_{\tilde{\D}} = \flip(\WD)$, without specifying $\s$, if
the precise choice of $\s$ does not matter.
\end{remark}

We can now describe the basic setup for implementing group operations on
the Jacobian, or more precisely on the classes of $\fieldk$-rational
divisors.  We will describe our algorithms in the context of the ``large
model'' of~\cite{KKM}, as well as a slight variant.  It is possible to
generalize our ideas to the ``medium'' and ``small'' models described in
that article, but the large model is sufficient to demonstrate the
asymptotic speedup of our new algorithms.

\begin{definition}
\label{definition3r.12}
The \emph{large model} of the curve $\C$ is defined as follows.  We
implicitly assume that $\g \geq 2$, although everything works (possibly
with some increase in degrees of divisors) for $\g \leq 1$.
\begin{enumerate}
\item
We choose a degree $\degd \geq 2\g$, with $\degd = O(\g)$ nonetheless, and
we fix once and for all an effective $\fieldk$-rational divisor $\D_0$
with $\deg \D_0 = \degd$.
\item
We define our basic line bundle by $\linebundleL = \OC(3\D_0)$, and
represent the spaces $\V$ and $\Vprime$ as well as the multiplication
map~$\mul$ using either \repA{} or \repB.  Note that $\Delta = 3\degd$.
\item
Given an effective $\fieldk$-rational divisor~$\D$, we say that $\D$ is
\emph{small} if $\deg \D = \degd$, and \emph{large} if $\deg \D = 2\degd$.
\item
If $\D$ is a small divisor, then let $x_\D$ be the linear
equivalence class of $\D - \D_0$ in the Jacobian of $\C$.  Then we
represent the element $x_\D$ of the divisor by the space $\WD$.
Similarly, if $\D$ is a large divisor, then define $x_\D$ to be the
linear equivalence class of $\D - 2\D_0$, and let the space $\WD$ represent
$x_\D$ .  
\item 
We calculate and store ahead of time the spaces $\W_{\D_0}$ and
$\W_{2\D_0}$, as well as \anigs{} for each space, and a specific $\s_0$,
unique up to a nonzero factor in $\fieldk$, such that $(\s_0)_\linebundleL
= 3\D_0$.  (Thus $\s_0$ corresponds to the element $1\in \fieldk(\C)$,
viewed as an element of $\HoOC{3\D_0}$.)
\item
If we need to perform the membership test of
Proposition/Algorithm~\ref{propalg4r.14}, or inflation as in 
Proposition/Algorithm~\ref{propalg3r.9}, then compute and store ahead of
time \anigs{} $\defl(\V)$ for $\V$ as mentioned above.
\end{enumerate}
\end{definition}

\begin{remark}
\label{remark3r.13}
Some assorted remarks:
\begin{enumerate}
\item
If the divisor $\D$ is small, then $\WD$ (respectively, $\WprimeD$)
has codimension $\degd$ in $\V$ (respectively, in $\Vprime$).
If $\D$ is large, then the codimension is $2\degd$.
Moreover, if $\D$ is small, then its complementary divisor
$\tilde{\D} = \flip(\D)$ is large, and vice-versa.  We see that $\D$ and
$\tilde{\D}$ represent inverse points on the Jacobian, since
$\D + \tilde{\D}$ is linearly equivalent to $3\D_0$. 
\item
We do not specifically need the spaces $\W_{\D_0}$ and $\W_{2\D_0}$.  We
can use instead spaces $\W_{\E_0}$ and $\W_{\E_1}$, where the divisor $\E_0$
is linearly equivalent to $\D_0$, and $\E_1 = \flip(\W_{\E_0}, \s_0)$ for
some nonzero choice of $\s_0 \in \W_{\E_0}$.  (It follows that $\E_1$ is
linearly equivalent to $2\D_0$.)
\item
When choosing the divisor $\D_0$ and the degree $\degd$, it is best to make
$\degd$ as small as possible, i.e., $\degd=2\g$ or perhaps $\degd = 2\g+2$
(which is useful in some contexts).  It may however be difficult in
practice to find effective divisors of a specific degree that are rational
over the base field $\fieldk$, especially if $\fieldk$ is a number field
(unless the curve $\C$ comes equipped with a known rational point).
\item
Assume that we start with a different representation of $\C$ before our
precomputation (e.g., as an equation for a singular plane curve birational
to $\C$).  We should also extend the precomputations of
Remark~\ref{remark2r.3} to compute some spaces $\WD$, for divisors
$\D$ that are supplied to us along with $\C$ (e.g., as formal sums of
points on the plane curve), and with which we wish to later do computations
in the Jacobian of $\C$.
\item
A side note: the divisor $\D_1$ in the definition of \repBzero{} and
Remark~\ref{remark2r.3} is $\D_1 = 3\D_0$.
\end{enumerate}
\end{remark}

We postpone until Section~\ref{section4r} a discussion of how to
quickly test whether a given subspace $\W \subset \V$, having the correct
dimension, actually is of the form $\WD$ for a small or large $\D$ ---
that membership test requires slightly different techniques from the other
algorithms, which in any case will be used much more often.  Instead, we
begin with a test for equality on the Jacobian.  Observe in this and our
later algorithms that we always perform a division by a deflation of a
subspace, i.e., using a small \igs{} instead of the entire subspace
representing a divisor.

\begin{propalg}[Equality of divisor classes]
\label{propalg3r.14}
Assume given two spaces $\WD$ and $\WE$, corresponding to divisors
$\D$ and $\E$ that are either both small or both large.  The the following
is a fast Las Vegas algorithm to test whether $\D$ and $\E$ are linearly
equivalent, i.e., whether $x_\D = x_\E$ on the Jacobian of $\C$:
\begin{enumerate}
\item
Take any nonzero $\s \in \WD$ and calculate
$\W = (\s \cdot \WE) \div \defl(\WD)$.
\item
Then $\D$ and $\E$ are linearly equivalent if and only if the
space $\W$ is nonzero.
\end{enumerate}
\end{propalg}
\begin{proof}
This is Theorem/Algorithm~4.1 of \cite{KKM}.  In brief, write
$(\s)_\linebundleL = \D + \tilde{\D}$, with $\D + \tilde{\D}$ linearly
equivalent to $3\D_0$.  Then
$\s \cdot \WE = \Wprime_{\D + \tilde{\D} + \E}$, so we obtain 
$\W = \W_{\tilde{\D} + \E}$ upon division.  Since
$\deg(\tilde{\D} + \E) = 3 \degd = \Delta$, the space
$\W_{\tilde{\D} + \E}$ is nonzero precisely when $\tilde{\D} + \E$ is 
linearly equivalent to $3\D_0$, which is equivalent to $\D$ and $\E$ being
linearly equivalent.  Note that $\deg(\tilde{\D} + \E)$ is larger than our
usual degree bounds; our computation of the space $\W_{\tilde{\D} + \E}$ is
nonetheless correct, as explained in~\cite{KKM}.
\end{proof}

For implementing group operations on the Jacobian, we shall be content with
describing one operation, ``addflip'':

\begin{definition}
\label{definition3r.14.5}
Given two elements $x, y$ in the Jacobian of $\C$ (actually, in any abelian
group that is written additively), we define their \emph{addflip} to be 
\begin{equation}
\label{equation3r.13}
\addflip(x,y) = -(x+y).
\end{equation}
Note that given this operation, it is of course immediate to compute
inverses, via $-x = \addflip(x,0)$, and hence to compute sums, via
$x+y = -\addflip(x,y)$.
\end{definition}

In the original large model from \cite{KKM}, we represented an element
of the Jacobian using only $\WD$ for a small divisor $\D$.  In that
context, we can implement the addflip as follows.

\begin{propalg}[Addflip of small divisors]
\label{propalg3r.15}
Assume given two subspaces $\WD$ and $\WE$, representing small divisors
$\D$ and $\E$, and elements $x_\D, x_\E$ of the Jacobian of $\C$.
divisors.  Then the following is a fast Las Vegas algorithm to compute a
space $\W_{\divF}$, for a suitable small divisor $\divF$, such that
$x_\divF = \addflip(x_\D,x_\E)$:
\begin{enumerate}
\item
Choose a nonzero $\s \in \WD$, and compute $\W_{\tilde{\D}} =
\flip(\WD,\s)$.  (Note that $\tilde{\D}$ is a large divisor.) 
\item
Compute $\W_{\D  + \E} = (\s \cdot \WE) \div \defl(\W_{\tilde{\D}})$.
(Note that $\D + \E$ is a large divisor.)
\item
Flip the result to obtain $\W_\divF = \flip(\W_{\D + \E})$.
\end{enumerate}
\end{propalg}
\begin{proof}
This is Proposition/Algorithm~4.3 of \cite{KKM}, using the second method of
adding divisors (Theorem/Algorithm~3.13 of that earlier article).  As in
Proposition/Algorithm~\ref{propalg3r.14} above, we have
$\s \cdot \WE = \Wprime_{\D + \tilde{\D} + \E}$, so our computation of
$\W_{\D+\E}$ is correct.  Step~3 shows that $\D + \E + \divF$ is linearly
equivalent to $3\D_0$, and hence $x_\D + x_\E + x_\divF = 0$ on the
Jacobian.
\end{proof}

\begin{remark}
\label{remark3r.15.5}
To evaluate $\addflip(0, x_\E)$, we of take $\D = \D_0$ and $\s = \s_0$.
This allows us to skip step~1, and simplify step~2, since we already know
a deflation of the space $\W_{\tilde{\D}} = \W_{2\D_0}$.
\end{remark}

As a variant, we can represent all elements on the Jacobian using large
divisors.  The resulting algorithm for addflip is given below.
Since $\D$ is now large, the space $\WD$ has smaller dimension than in our
original large model.  This will make some computations faster, especially
since we do fewer basic operations in this algorithm than in
Proposition/Algorithm~\ref{propalg3r.15}.

\begin{propalg}[Addflip of large divisors]
\label{propalg3r.16}
Given two elements $x_\D,x_\E$ of the Jacobian, represented by $\WD, \WE$
for large divisors $\D, E$, we can compute $\W_{\divF}$ for a large divisor
$\divF$ that represents $x_\divF = \addflip(x_\D,x_\E)$ by the following
fast Las Vegas algorithm:
\begin{enumerate}
\item
Compute $\W_{\tilde{\D}} = \flip(\WD)$. 
(Note that $\tilde{\D}$ is a small divisor.)
\item
Choose a nonzero $\s \in \WE$, so $(\s) = \E + \tilde{\E}$.
Compute $\W_{\tilde{\D} + \tilde{\E}}
    = (\s \cdot \W_{\tilde{\D}}) \div \defl(\WE)$. 
\item
Our desired result is $\W_{\divF} = \W_{\tilde{\D} + \tilde{\E}}$.
\end{enumerate}
\end{propalg}
\begin{proof}
The inverses $-x_\D$ and $-x_\E$ in the Jacobian are given by the linear
equivalence classes of $\tilde{\D} - \D_0$ and $\tilde{\E} - \D_0$.  Thus
the divisor $\divF = \tilde{\D} + \tilde{\E}$ represents $-x_\D - x_\E$.
\end{proof}

\section{Randomly selecting \anigs, with verification; membership test}
\label{section4r}

In the first part of this section, we are given an effective
$\fieldk$-rational divisor $\D$ for which $\WD$ is base point free, and we
let $\W \subset \WD$ be a subspace whose divisor of common zeros is $\D$
(in most applications, $\W = \WD$).  We wish to study the probability that
a suitable random selection of $\s_1, \dots, \s_\h \in \W$ is \anigs{} for
$\D$.  In order to clarify what is going on, we shall work with the line
bundle $\linebundleM = \linebundleL(-\D)$.  Then we can view $\W$ as a
base point free subspace of $\HoM$, more precisely as a base point free
linear series of the line bundle $\linebundleM$.  We hence wish to
determine the probability that there is no point common to all the divisors
$(\s_1)_\linebundleM, \dots, (\s_\h)_\linebundleM$.

\begin{lemma}
\label{lemma4r.1}
Let $\linebundleM$ be a base point free line bundle on $\C$.  Let $\W
\subset \HoM$ be a base point free subspace.  Fix a nonzero element
$\s_1 \in \W$. 
\begin{enumerate}
\item
There exist proper subspaces $H_1, \dots, H_\ell \subsetneq \W$, with
$\ell \leq \deg \linebundleM$, with the following property:
\begin{equation}
\label{equation4r.1}
\{\s_2 \in \W \mid \{\s_1,\s_2\} \text{ is NOT \anigs{} for } \HoM \}
     = H_1 \union \dots \union H_\ell.
\end{equation}
\item
More generally, let $\h \geq 2$, and view a selection of $\s_2, \dots,
\s_\h \in \W$ as a tuple $(\s_2, \dots, \s_\h)$ in the vector space
$\W^{\h-1}$.  Then, with the same $\{H_i\}$ as in part~1,
\begin{equation}
\label{equation4r.2}
\begin{split}
&\{(\s_2, \dots, \s_\h) \in \W^{\h-1} \mid \{\s_1,\dots,\s_\h\}
     \text{ is NOT \anigs{} for } \HoM \} \\
&\qquad = (H_1)^{\h-1} \union \dots \union (H_\ell)^{\h-1}.\\
\end{split}
\end{equation}
\end{enumerate}
\end{lemma}
\begin{proof}
Let $\pointP_1, \dots, \pointP_\ell \in \C(\fieldkbar)$ be the
\emph{distinct} points where $\s$ vanishes.  Thus
$\ell \leq \deg\linebundleM$.  Define $H_i$ to be the $\fieldk$-rational
subspace $\{\tu \in \W | v_{\linebundleM,\pointP_i}(\tu) \geq 1\}$
of sections vanishing at $\pointP_i$.  Since $\W$ is base point
free, we have $H_i \subsetneq \W$.  Then both sides of~\eqref{equation4r.2}
express the fact that all of $\s_2, \dots, \s_\h$ also vanish at one of the
$\pointP_i$. 
\end{proof}

The next lemma is an abstract statement about linear algebra; we have
adapted it from a result in~\cite{BrawleyGao}.
\begin{lemma}
\label{lemma4r.2}
Let $\W$ be a vector space over $\fieldk$, with basis
$\{w_1, \dots, w_r\}$.  Take a finite subset $\Sigma \subset \fieldk$, and
consider $\Sigma$-random elements of $\W$ in the sense of
Definition~\ref{definition2r.9}.  Let $H_1, \dots, H_\ell \subsetneq \W$ be
proper subspaces.
\begin{enumerate}
\item
For a $\Sigma$-random element $\tu \in \W$,
\begin{equation}
\label{equation4r.3}
\Pr(\tu \in H_1 \union \dots \union H_\ell) \leq \ell/\abs{\Sigma}.
\end{equation}
\item
For a tuple $(\tu_1, \dots, \tu_j) \in \W^j$ of independent
$\Sigma$-random elements $\tu_1, \dots, \tu_j \in \W$,  
\begin{equation}
\label{equation4r.4}
\Pr\Bigl((\tu_1, \dots, \tu_j) \in
    (H_1)^j \union \dots \union (H_\ell)^j\Bigr)
 \leq \ell/\abs{\Sigma}^j. 
\end{equation}
\end{enumerate}
\end{lemma}
\begin{proof}
Both statements easily reduce to the case $\ell=1$, so we assume from now
on that we only have one subspace $H = H_1 \subsetneq\W$.
We can find an $(r-1)$-dimensional hyperplane $H' \subset \W$ containing
$H$.  Hence there exist constants $a_1, \dots, a_r \in \fieldk$, not all
zero, such that 
\begin{equation}
\label{equation4r.5}
\tu = c_1 w_1 + \dots + c_r w_r \in H \implies
\tu \in H' \iff
   a_1 c_1 + \dots + a_r c_r = 0.
\end{equation}
Without loss of generality, say that $a_1 \neq 0$.  Then for every choice
of values of $c_2, \dots, c_r \in \Sigma$, there exists exactly one value
of $c_1 \in \fieldk$ for which $\tu \in H'$, hence at most one value 
of $c_1$ for which $\tu \in H$; it is furthermore possible that
$c_1 \notin \Sigma$.  So at most $\abs{\Sigma}^{r-1}$ choices of tuples
$(c_1, \dots, c_r) \in \Sigma^r$ lead to $\tu \in H$, whence
$\Pr(\tu\in H) \leq 1/\abs{\Sigma}$.  It follows that
$\Pr\Bigl((\tu_1, \dots, \tu_j) \in H^j\Bigr) \leq 1/\abs{\Sigma}^j$.  This
proves our result.
\end{proof}

Combining the above two lemmas, we immediately obtain:

\begin{proposition}
\label{proposition4r.3}
Keep the assumptions and notation of Lemmas \ref{lemma4r.1}
and~\ref{lemma4r.2} above.  Take $0 < \eta < 1$, and define
\begin{equation}
\label{equation4r.6}
\h = 1 + \left\lceil 
            (\log \deg \linebundleM - \log \eta)/ \log \abs{\Sigma}
      \right\rceil.
\end{equation}
For a fixed nonzero $\s_1 \in \W$, let $\s_2, \dots, \s_\h \in \W$ be
independently chosen $\Sigma$-random elements.  Then
\begin{equation}
\label{equation4r.7}
\Pr\Bigl(\{\s_1, \dots, \s_\h\} \text{ is \anigs{} for } \HoM \Bigr) 
  \geq
 1 - \eta.
\end{equation}
\end{proposition}
\begin{proof}
Immediate, once we note that $j = \h-1$ in our previous notation, and that
$\ell \leq \deg \linebundleM$.
\end{proof}

\begin{corollary}
\label{corollary4r.4}
If $\fieldk$ is infinite, then every base point free subspace $\W$ contains
\anigs{} with two elements.
\end{corollary}

We are now ready for a more precise statement about random sections
giving \anigs, when $\fieldk$ is a finite field.  We thus take $\Sigma =
\fieldk$; a $\Sigma$-random element of a vector space $\W$ is thus a
random element of the finite set $\W$, chosen using the uniform
distribution.  We first note two simple facts.

\begin{lemma}
\label{lemma4r.6}
Assume that $\fieldk = \Fq$.
For $\ell \geq 1$, let $N_\ell$ be the number of degree $\ell$
irreducible divisors on $\C$ (i.e., divisors of the form $\D = \pointP_1 +
\dots + \pointP_\ell$, where the $\ell$ points $\{\pointP_1, \dots,
\pointP_\ell\}$ are a single Galois orbit).  Then
\begin{equation}
\label{equation4r.12}
N_\ell \leq \frac{1}{\ell} (\q^\ell + 1 + 2\g \q^{\ell/2}).
\end{equation}
\end{lemma}

\begin{proof}
The $N_\ell$ irreducible divisors give rise to $\ell N_\ell$ distinct
$\Fql$-rational points on $\C$.  
However, $\abs{C(\Fql)} \leq \q^\ell +  1 + 2\g \q^{\ell/2}$ by the simplest
form of the Weil bounds (see for example Appendix~C of~\cite{Hartshorne}).
\end{proof}

\begin{lemma}
\label{lemma4r.7}
Assume that $\fieldk = \Fq$, and that $\deg \linebundleM = \dbar + 2\g - 1$
with $\dbar \geq 1$.  Choose random $\s_1, \dots, \s_\h \in \HoM$
independently with the uniform distribution.  Then the probability that the
sections have a common zero (i.e., that they are \emph{not} \anigs) is at
most 
\begin{equation}
\label{equation4r.13}
N_1 \q^{-\h} + N_2 \q^{-2\h} + \dots + N_{\dbar} \q^{-\dbar\h} +
N_{\dbar+1} \q^{-\dbar\h} + \dots + N_{\dbar + 2\g - 1} \q^{-\dbar\h}.
\end{equation}
\end{lemma}

\begin{proof}
For each irreducible divisor $\D$, the probability that a given section
vanishes at $\D$ is $\abs{\HoMm{\D}} / \abs{\HoM} = \q^{-c}$, where $c$ is
the codimension of $\HoMm{\D}$ in $\HoM$.  Thus the probability that $\h$
sections all vanish at $\D$ is $\q^{-ch}$.  Now by Riemann-Roch,
we have that $c = \deg \D$ when $1 \leq \deg \D \leq \dbar$, and
$c \geq \dbar$ when $\deg \D \geq \dbar$.  Moreover, we know that if
$\deg \D \geq \dbar + 2\g$, then $\HoMm{\D} = \{0\}$, so in that case
simultaneous vanishing 
at $\D$ can happen only if all the sections are identically zero ---
but we have already accounted for this situation in considering
divisors of smaller degree.  Adding up for all irreducible $\D$ the
probability that the sections simultaneously vanish at $\D$ yields the
upper bound~\eqref{equation4r.13}.
\end{proof}

\begin{remark}
\label{remark4r.8}
In the above proof, we have not tried to bound the ``overcounting''; for
example, if $\D_1$ and $\D_2$ are distinct irreducible divisors, then
we have counted twice the contribution to \eqref{equation4r.13} of the
probability that the sections all vanish at $\D_1 + \D_2$.
Heuristically, at least when $\dbar \to \infty$, the events of vanishing
at two (or more) divisors $\D_1, \D_2$ should be independent, with
probabilities $\q^{-\deg \D_1}$ and $\q^{-\deg \D_2}$.  So for $\dbar$
large, a good heuristic estimate of the probability that $\h$ random
sections do not yield \anigs~is given by
\begin{equation}
\label{equation4r.14}
1 - \prod_{\ell = 1}^\infty ( 1 - \q^{-\ell \h} )^{N_\ell}
 = 1 - \frac{1}{Z_\C(\h)},
\end{equation}
where $Z_\C(s)$ is the zeta function of $\C$.  This is analogous to a
standard elementary statement that the ``probability'' that two
integers $m,n \in \Z$ are relatively prime (i.e., that $\{m,n\}$ is
\anigs!) is $\prod_{p \text{ prime}} (1-p^{-2}) = 1/\zeta(2) =
6/\pi^2$.  Now $Z_\C(s)$ is a rational function of $\q^{-s}$, and its
expansion near $\q^{-s}=0$ (i.e., as $s \to \infty$) gives us
$1- 1/Z_\C(\h) = N_1 \q^{-\h} + O(\q^{-2\h})$.  Thus if we want this
quantity to be less than $\eta$,  we can try the heuristic approximation
$\h = \lceil \log(N_1/\eta)/\log \q) \rceil$.  Now
$N_1 \leq \q + 1 + 2\g\sqrt{\q}$, so if we fix $\q$ and let $\g$ become
large, we obtain a value $\h \approx \log(2\g\sqrt{\q}/\eta)/\log(\q)
= O(1 + \log(\g/\eta)/\log \q)$, in line with our results.
\end{remark}

We can now state and prove our result Proposition~\ref{proposition4r.9} for
finite fields.  Even though 
our algorithms rely on the simpler Proposition~\ref{proposition4r.3}, the
significance of Proposition~\ref{proposition4r.9} is that the value of $\h$
given below does not depend on $\dbar$, once $\dbar$ is comparable to or
larger than $\g$.  Also note that if $\g$ or $\q$ is large, then the
constant $6$ in~\eqref{equation4r.15} can be reduced significantly.
However, the result of Proposition~\ref{proposition4r.9} only works if we
randomly select our sections from the entire space $\HoM$, and not a
subspace $\W$.

\begin{proposition}
\label{proposition4r.9}
Assume that $\fieldk = \Fq$, and that $\deg \linebundleM = \dbar + 2\g - 1$
with $\g \geq 1$ and $\dbar \geq 2$.  Let $0 < \eta < 1$, and define
\begin{equation}
\label{equation4r.15}
\h = \max \left( 1 + \left\lceil \frac{2\g-1}{\dbar-1} \right\rceil,
           1 + \left\lceil
                    \frac{\log (6\g/\eta)}{\log \q} 
               \right\rceil
      \right).
\end{equation}
Then a uniform random choice of $\h$ sections from $\HoM$ is \anigs~with
probability $> 1-\eta$.
\end{proposition}

\begin{proof}
By Lemmas \ref{lemma4r.6} and~\ref{lemma4r.7}, the probability of
\emph{not} being \anigs~is bounded above by the quantity
\begin{equation}
\label{equation4r.16}
P = \sum_{\ell = 1}^{\dbar}
         \frac{1}{\ell} (\q^\ell + 1 + 2\g \q^{\ell/2}) \q^{-\ell \h}
  + \sum_{\ell = \dbar+1}^{\dbar+2\g-1}
         \frac{1}{\ell} (\q^\ell + 1 + 2\g \q^{\ell/2}) \q^{-\dbar \h}.
\end{equation}
We wish to show that $P < \eta$.  We use the following elementary
estimates that hold for $N \geq M \geq 1$, $\q \geq 2$, and $\sigma
\geq 1$:
\begin{equation}
\label{equation4r.17}
\begin{split}
\sum_{\ell=1}^M \frac{1}{\ell} \q^{-\ell \sigma}
          < \q^{-\sigma} + \frac{\q^{-2\sigma}}{2(1-\q^{-\sigma})}
          \leq 1.5 \q^{-\sigma};
\qquad &
\sum_{\ell=M}^N \frac{1}{\ell} \leq \frac{N-M+1}{M};
 \\
\sum_{\ell=M}^N \frac{1}{\ell} \q^\ell
      \leq \frac{1}{M} \cdot \frac{\q^N}{1-\q^{-1}}
      \leq \frac{2\q^N}{M};
\qquad &
\sum_{\ell=M}^N \frac{1}{\ell} \q^{\ell/2} <
      \frac{3.5 \q^{N/2}}{M}.
 \\
\end{split}
\end{equation}
(The constant $3.5$ is a simple upper bound for $1/(1-\q^{-1/2})$ when
$\q \geq 2$.)
From these, we easily estimate that
\begin{equation}
\label{equation4r.18}
\begin{split}
 P  & < 1.5 \left[ \q^{-(\h-1)} + \q^{-\h} + 2\g \q^{-(\h-1/2)}
             \right] \\
    & \qquad +
     \frac{\q^{-\dbar\h + \dbar + 2\g - 1}}{\dbar+1}
     \left[ 2 + (2\g - 1)\q^{-(\dbar+2\g-1)} + 7\g \q^{-(\dbar+2\g-1)/2}
     \right]. 
\end{split}
\end{equation}
(Note that $\h - 1 \geq 1$.)  Equation~\eqref{equation4r.15}
now implies that $\q^{1-\h} \leq \eta/6\g$, and also that
$ - \dbar\h + \dbar + 2\g - 1 \leq 1 - \h < 0$.  Since furthermore
$\dbar + 2\g - 1 \geq 3$ , we obtain
\begin{equation}
\label{equation4r.19}
\begin{split}
P & <
      1.5 \left[
 \frac{\eta}{6\g} + \frac{\eta}{6\g\q} + \frac{2\eta}{6\sqrt{\q}}
        \right]
 +  \frac{\eta}{6(\dbar+1)\g} \left[
               2 + \frac{2\g - 1}{\q^3} + \frac{7\g}{\q^{1.5}}
                          \right]
 \\
& \leq
 \frac{\eta}{6}
    \left[
       \frac{1.5}{\g} + \frac{1.5}{\g\q} + \frac{3}{\sqrt{\q}}
        + \frac{ (2-\q^{-3})\g^{-1} +  2\q^{-3} + 7\q^{-1.5}}
               {\dbar+1}
    \right]
  < \eta,   
\end{split}
\end{equation}
since $\q \geq 2$, $\g \geq 1$, and $\dbar \geq 2$.  This gives the
desired result.
\end{proof}

Our second topic in this section is to discuss how to verify whether our
random selection of sections is indeed \anigs.  The same techniques also
give our algorithm for membership testing.  We prove both these results
after two preliminary lemmas.  We return to considering a line bundle
$\linebundleL$, of degree $\Delta \geq 2\g + 2$, and subspaces of the form 
$\WD \subset \V, \WprimeD \subset \Vprime$ for effective $\fieldk$-rational
divisors $\D$.

\begin{lemma}
\label{lemma4r.11}
Let $\D$ be an effective divisor for which $\WD$ is base point free.  Then
$\deg \D \leq \Delta$.  Morover, we have the following 
relation between $\deg \D$ and the codimension of $\WD$ in $\V$:
\begin{enumerate}
\item
If $\codim \WD \leq \Delta - 2\g$, then $\deg \D = \codim \WD$.
\item
If $\codim \WD \geq \Delta - 2\g + 1$, then
$\deg \D \geq \Delta - 2\g + 1$.
\end{enumerate}
\end{lemma}
\begin{proof}
The first statement follows because $\D$ is a ``factor'' of the divisor of
any nonzero $\s \in \WD$, but $\deg (\s)_{\linebundleL} = \Delta$.  The
statements about the codimension are straightforward (extend scalars to
$\fieldkbar$, start with $\D=0$, and add one point at a time to $\D$).
\end{proof}

\begin{lemma}
\label{lemma4r.12}
Assume given nonzero $\tu_1, \tu_2 \in \V$ such that $\D$ is the divisor
of common zeros of $\tu_1, \tu_2$.  Define
$\Wprime = \tu_1 \cdot \V + \tu_2 \cdot \V$.  Then $\Wprime \subset
\WprimeD$, and the codimension of $\Wprime$ in $\Vprime$ satisfies
\begin{equation}
\codim \Wprime = \dim \HoOC{\D} - 1 + \g
               = \deg \D + \dim H^1\bigl(\OC(\D)\bigr).
\label{equation4r.20}
\end{equation}
In particular, if $\deg \D \geq 2\g-1$, then $\Wprime = \WD$.
\end{lemma}
\begin{proof}
Write $(\tu_1)_\linebundleL = \D + \E_1$ and $(\tu_2)_\linebundleL 
= \D + \E_2$, where $\E_1$ and $\E_2$ are disjoint effective divisors.  Now
$\tu_1 \cdot \V = \Wprime_{\D + \E_1}$ and 
$\tu_2 \cdot \V = \Wprime_{\D + \E_2}$; hence trivially
$\Wprime \subset \WprimeD$. 
We now use
$\codim \Wprime = \codim (\tu_1 \cdot \V) + \codim (\tu_2 \cdot V) - \codim
(\tu_1\cdot\V \intersect \tu_2\cdot\V)$
to show~\eqref{equation4r.20}.
By construction, 
$\tu_1\cdot\V \intersect \tu_2\cdot\V = \Wprime_{\D + \E_1 + \E_2}$.
Now $\D + \E_1$ and $\D + \E_2$ are in the linear equivalence class of
$\linebundleL$, so 
$\linebundleL^{\tensor2}(-\D - \E_1 - \E_2) \isomorphic \OC(\D)$.
Therefore
$\dim \Wprime_{\D + \E_1 + \E_2} = \dim \HoOC{\D}$, and its codimension is
$\delta' - \dim \HoOC{\D} = 2\Delta + 1 - \g - \dim \HoOC{\D}$.  On the
other hand, both $\tu_1 \cdot \V$ and $\tu_2 \cdot \V$ have codimension
$\Delta$ in $\Vprime$.  This proves~\eqref{equation4r.20}.  As for the
last statement, note that the assumption on $\deg \D$ implies that $\codim
\Wprime = \deg \D$.  However, we always have $\Wprime \subset \WprimeD$,
and moreover $\codim \WprimeD = \deg \D$ (use Lemma~\ref{lemma4r.11}
to get $\deg \D \leq \Delta \leq 2\Delta - 2\g$).
Thus $\Wprime = \WD$, as desired.
\end{proof}

\begin{proposition}
\label{proposition4r.13}
The criterion of part~2 of Theorem~\ref{theorem2r.10} is correct.
\end{proposition}
\begin{proof}
It is enough to prove the statement after extending scalars to
$\fieldkbar$ (any infinite field will do).  Let $S \subset \WD$ be the
subspace spanned by $\s_1, \dots, \s_\h$.  Then the divisor of common
zeros of $S$ is $\D + \divF$ for some effective divisor $\divF$, and
$\{\s_1, \dots, \s_\h\}$ is \anigs{} for $\D$ if and only if $\divF = 0$.
By Corollary~\ref{corollary4r.4}, there exist $\tu_1, \tu_2 \in S$ whose
divisor of common zeros is also $\D + \divF$.  We have
\begin{equation}
\label{equation4r.21}
 \tu_1 \cdot \V + \tu_2 \cdot \V \subset
 \s_1 \cdot \V + \dots + \s_\h \cdot \V \subset
 \Wprime_{\D+\divF} \subset
 \WprimeD.
\end{equation}
Apply the final statement of Lemma~\ref{lemma4r.12} to the divisor
$\D+\divF$, whose degree is at least $2\g - 1$ by the assumption on
$\deg \D$; we conclude that $\s_1 \cdot \V + \dots \s_\h \cdot \V =
\Wprime_{\D + \divF}$, with codimension $\deg \D + \deg \divF$.  This
yields the desired result.
\end{proof}

\begin{propalg}[Membership test]
\label{propalg4r.14}
Given a subspace $\W \subset \V$, write $c = \codim \W$ (in $\V$), and
assume that $2\g \leq c \leq \deg \linebundleL - 2\g$.  Define
$\h = 1 + \lceil \log 2\Delta / \log \abs{\Sigma}\rceil$.
Let $\D$ be the divisor of common zeros of $\W$, so $\W \subset \WD$.  Then
the following is a fast algorithm to check if $\W = \WD$, under the
assumption that we have precomputed \anigs{} $\defl(\V)$ for $\V$ as in
Lemma/Algorithm~\ref{lemmalg3r.8}: 
\begin{enumerate}
\item
Select $\s_1, \dots, \s_\h \in \W$ in the usual way (take any 
$\s_1 \neq 0$, and choose the rest $\Sigma$-randomly), and calculate
\begin{equation}
\label{equation4r.22}
U' = \s_1 \cdot \V + \dots + \s_\h \cdot \V.
\end{equation}
Write $c' = \codim U'$ (in $\Vprime$).  If $c' > c$,  
then go back to step~1.  Else, if $c' < c$, then conclude that 
$\W \neq \WD$ and stop.  Otherwise (if $c'=c$), continue.
\item
Compute $U = U' \div \defl(\V)$.  If $U = \W$, then conclude
that $\W = \WD$.  Otherwise, conclude that $\W \neq \WD$.
\end{enumerate}
\end{propalg}
\begin{proof}
By statement~(1) of Lemma~\ref{lemma4r.11}, we have
$c \geq \codim \WD = \deg D$.  Our choice of $\h$ (which is still
$O(\g^\eps)$) implies that $\{\s_1, \dots, \s_\h\}$ is \anigs{} for
$\D$ with probability at least $1/2$, independently of $\deg \D$.
As in Proposition~\ref{proposition4r.13}, we extend scalars to
$\fieldkbar$, and write $S = \text{span}\{\s_1, \dots, \s_\h\}$, with
divisor of common zeros $\D + \divF$; we have $\divF = 0$ at least half
the time.  Again let $\tu_1, \tu_2 \in S$ have divisor of common zeros
$\D + \divF$, to obtain the same inclusions as in~\eqref{equation4r.21}.

We now discuss what happens in the two cases $\W = \WD$ and $\W \neq \WD$:
\begin{enumerate}
\item
If $\W = \WD$, then $\deg \D = c$, and we obtain as in
Proposition~\ref{proposition4r.13} that $c' = \deg (\D + \divF)$; thus
$c' - c = \deg \divF \geq 0$.  We therefore repeat the loop in step~1
at most twice on average until we have $\divF = 0$, at which point we also
obtain $U' = \WprimeD$.  It follows that the division in step~2
computes $U = \WD$, so the test correctly concludes that $\W = \WD$. 
\item
If $\W \neq \WD$, then $\deg \D < c$.  We distinguish four scenarios: 
\begin{enumerate}
\item 
$\deg (\D + \divF) \leq 2\g - 1$:
Let $\divF'$ be any effective
divisor for which $\deg (\D + \divF + \divF') = 2\g - 1$, and note that
$\dim \HoOC{\D+\divF} \leq \dim \HoOC{\D+\divF+\divF'} = \g$.  By
Lemma~\ref{lemma4r.12}, we know that 
$c' \leq \codim{\tu_1\cdot \V + \tu_2 \cdot \V} \leq 2\g-1 < c$.  
Hence the test correctly concludes that $\W \neq \WD$ in step~1.
\item 
$2\g - 1 \leq \deg (\D + \divF) < c$:
in this and in the following scenarios, we have $c' = \deg(\D + \divF)$ and
$U' = \Wprime_{\D + \divF}$, as in the previous proposition.  So in this
particular scenario, $c' < c$, and we conclude that $\W \neq \WD$ in step~1.
\item 
$\deg (\D + \divF) = c$:
This time, $c' = c$, and we move on to step~2, where we compute
$U = \W_{\D + \divF}$.  It follows that $\W \not\subset U$, because $\D$ is
the divisor of common zeros of $\W$, whereas $\divF \neq 0$ (its degree is
$c - \deg D$).  Thus step~2 correctly concludes that $\W \neq \WD$.
\item 
$c < \deg (\D + \divF)$:
Here $c' > c$, so we repeat step~1.  This happens less
than half the time, since if $\divF = 0$ we must have already landed in
scenario a or~b above.  Thus we loop in step~1 at most twice on average.
\end{enumerate}
\end{enumerate}
\end{proof}

\begin{remark}
\label{remark4r.15}
Our original ``slow'' algorithm for testing whether $\W = \WD$,
Theorem/Algorithm 3.14 of~\cite{KKM}, was to compute $\flip(\W)$ 
and to see if the result had the expected dimension.  There, the flip
was implemented using a division by a basis for $\W$, which was \anigs{}
for $\D$.  We unfortunately cannot do the same using a random selection of
$\h = O(\g^\eps)$ elements from $\W$ as our \igs, because we would not be
able to quickly verify whether our random selection actually was \anigs{}
(we do not know $\deg \D$ in advance, and it is moreover likely that the
$\deg \D < 2\g - 1$).
\end{remark}

\section{Converting from \repA{} to \repB}
\label{section5r}

Our goal in this section is to give a brief sketch, under some conditions
on $\fieldk$ given below, of how we can convert a curve $\C$ given using
\repA{} into a description of $\C$ using \repB.  This is a precomputation
that we only need to do once, so we will be satisfied with an efficient
algorithm (as defined below), which is essentially polynomial time, but not
necessarily of complexity $O(\g^{3+\eps})$.

We emphasize, however, that if it is at all possible to find enough
points in $\C(\fieldk)$ so as to use the simpler form \repBzero, then we
should do so, even if we do not bother with fast linear algebra.  For
example, this should not pose a problem if $\fieldk = \Fq$ with $\q$ very
large compared to $\g$, since then $\abs{\C(\Fq)}$ is comparable to $\q$.

\emph{In this section, we maintain the following two assumptions about our
field $\fieldk$.} 
Both of these assumptions hold if $\fieldk$ is a finite field or a number
field.
\begin{enumerate}
\item
The field $\fieldk$ is perfect.
\item
There exists an efficient algorithm to compute the primary decomposition
(including finding the radical) of a finite-dimensional $\fieldk$-algebra
$\mathcal{A}$.
\end{enumerate}
The second condition can, nontrivially, be replaced by our being able to
efficiently factor (univariate) polynomials in $\fieldk[x]$.  Here an
\emph{efficient algorithm} means that if $N = \dim \mathcal{A} = O(\g)$,
then we have a Las Vegas algorithm with an expected complexity that is
polynomial in $\g$, where we need to measure complexity in terms of
\emph{both} field operations and factorizations of degree $O(N)$
polynomials in $\fieldk[x]$.  As examples of algorithms for primary
decomposition and the computation of radicals, we mention the articles
\cite{EberlyGiesbrecht}, \cite{Kemper}, and~\cite{DeckerGreuelPfister}, and
the articles cited in their bibliographies.

For an extended treatment of the material in this section, including many
details omitted here as well as a fairly self-contained algorithm
for primary decomposition, the reader is referred to Sections 6 and~7 of
\cite{KKMextendedpreprint}.

Starting from \repA, we can as before produce the projective coordinate
ring of $\C$, as in Proposition~\ref{proposition2r.1} and
Lemma/Algorithm~\ref{lemmalg3r.8}: this is
\begin{equation}
\label{equation5r.1}
\Directsum_{n \geq 0} H^0 \bigl( \linebundleL^{\tensor n} \bigr)
\isomorphic
\fieldk[T_1, \dots, T_\delta] / I_\C.
\end{equation}
We choose the divisor $Z$ for \repB{} to be
\begin{equation}
\label{equation5r.2}
Z = 3 (T_1)_\linebundleL = (T_1^3)_{\linebundleL^{\tensor 3}}.
\end{equation}
Note that if we view $\linebundleL = \OC(\D_1)$, then we can take $Z =
3\D_1$.  Here $\deg Z = 3\Delta$, which allows us to faithfully represent
elements of $\V$ and $\Vprime$ by their ``values'' at $Z$.  The values in
question belong to the algebra $\mathcal{A} = H^0(\mathcal{O}_Z)$, which
has dimension $3\Delta = O(\g)$.

\begin{proposition}
\label{proposition5r.1}
We can efficiently find a description of $\mathcal{A}$ in terms of a
basis and a multiplication table for $\mathcal{A} \times \mathcal{A} \to
\mathcal{A}$ (similarly to \eqref{equation2r.6} and~\eqref{equation2r.7}).
In the process, we also obtain the images of $T_1, \dots, T_\delta$ as 
linear combinations of our basis for $\mathcal{A}$, which allows us to
identify $\V$ with a specific subspace of $\mathcal{A}$.
\end{proposition}
\begin{proof}[Sketch of proof]
View $Z$ as a zero-dimensional
subscheme of the projective space containing $\C$.  Its projective
coordinate ring is then
$\fieldk[T_1, \dots, T_\delta]/\bigl(I_\C + (T_1^3)\bigr)$.  
We however need to find the \emph{affine} coordinate ring $\mathcal{A}$ of
$Z$.  We first deal with an easy case, when $\{T_1, T_2\}$ is \anigs{} for
$\V$.  (This can be arranged, for example, if $\fieldk$ has at least
$2\Delta$ elements, since we can then choose $T_2$ randomly with a good
chance of getting \anigs, which we can verify as in
Lemma/Algorithm~\ref{lemmalg3r.8}.)  In this easy case, the scheme $Z$ lies
entirely in the affine open subset of projective space given by
$T_2 \neq 0$, so we can take
\begin{equation}
\label{equation5r.3}
\mathcal{A} 
  = H^0(\mathcal{O}_Z) 
  = k[T_1, \dots, T_\delta]/\bigl( I_\C + (T_1^3) + (T_2 - 1) \bigr).
\end{equation}
The images of $T_1, \dots, T_\delta$ in $\mathcal{A}$ are the obvious
ones.  We can find a basis and multiplication table for $\mathcal{A}$ using
Gr\"obner bases, or by a more direct approach that uses our linear algebra
algorithms on subspaces of $H^0 \bigl(\linebundleL^{\tensor n} \bigr)$ for
$n \leq 8$, which more clearly shows that the algorithm is efficient.

As for the more general case, we need to consider all affine open subsets
given by $T_j \neq 0$ for $2 \leq j \leq \delta$.  (It suffices in fact to
consider $2 \leq j \leq \h$, where $\{T_1, T_2, \dots, T_\h\}$ is \anigs{}
for $\V$.)  For each such $j$, we form the quotient ring
of~\eqref{equation5r.3}, but with $T_j$ instead of $T_2$.  The quotient
ring is then $H^0(\mathcal{O}_{Z_j})$, where $Z_j$ is the portion of $Z$
lying in the affine open set $\{T_j \neq 0\}$.  It is then possible to put
the $\{H^0(\mathcal{O}_{Z_j})\}$ together, while eliminating redundancy
from the intersections $Z_i \intersect Z_j$, to obtain $\mathcal{A}$ as
above.  (Roughly speaking, remove $Z_2$ from $Z$ using a division, then
remove any part of $Z_3$ from what remains, and so forth, finding the
affine algebra of each piece; then $\mathcal{A}$ is the product of these
partial affine algebras.)  All this can again be done using
only linear algebra on subspaces of 
$H^0 \bigl(\linebundleL^{\tensor n} \bigr)$ for $n \leq 8$.
\end{proof}

Now that we have represented $\mathcal{A}$ in a form suitable for
computation, we use our ability to find primary decompositions to decompose
$\mathcal{A}$ into a product of local Artinian $\fieldk$-algebras:
\begin{equation}
\label{equation5r.4}
\mathcal{A} = \mathcal{B}_1 \times \dots \times \mathcal{B}_r.
\end{equation}
This decomposition corresponds to writing $Z = e_1 Y_1 + \dots e_r Y_r$ for
distinct irreducible divisors $Y_i$ (cf.~Lemma~\ref{lemma4r.6}). Thus the
above decomposition expresses the canonical isomorphism
\begin{equation}
\label{equation5r.5}
H^0(\mathcal{O}_Z) \isomorphic
H^0(\mathcal{O}_{e_1 Y_1}) \times \dots \times H^0(\mathcal{O}_{e_r Y_r}).
\end{equation}
Let $R$ be the affine coordinate ring of any fixed open subset of $\C$ that
contains $Z$.  Then each irreducible divisor $Y_i$ corresponds to a maximal
ideal $P_i$ of the Dedekind domain $R$, and $Z$ corresponds to the ideal
$J = P_1^{e_1} \cdots P_r^{e_r}$.
Then the above decompositions are just the Chinese Remainder Theorem:
\begin{equation}
\label{equation5r.6}
R/J \isomorphic R/P_1^{e_1} \times \dots \times R/P_r^{e_r}.
\end{equation}
We will use the existence of $R$ and the $P_i$ to clarify our exposition,
but we point out that we do not compute $R$ at all; all our calculations
occur in the finite-dimensional algebra $\mathcal{A}$ and in certain vector
space subquotients such as the $\mathcal{B}_i$.

Specifically, the primary decomposition algorithm gives us an explicit
basis for each $\mathcal{B}_i$, viewing $\mathcal{B}_i$ as a
$\fieldk$-subspace of $\mathcal{A}$.  We simultaneously obtain, via the
computation of the radical, a basis for the maximal ideal $\mathfrak{p}_i$
of $\mathcal{B}_i$; here the inclusion
$\mathfrak{p}_i \subset \mathcal{B}_i$ corresponds to 
$P_i/P_i^{e_i} \subset R/P_i^{e_i}$.
Write $L_i$ for the residue field
$\mathcal{B}_i/\mathfrak{p}_i \isomorphic R/P_i$; in terms of $f_i =
[L_i:\fieldk]$, we have $\dim_\fieldk \mathcal{B}_i = e_i f_i$.   
Using our multiplication table for $\mathcal{A}$, we can easily implement
the ring operations in either $\mathcal{B}_i$ or $L_i$.  We can also
determine any $\fieldk$-linear dependencies between the elements of any
finite subset $\{\beta_1, \dots, \beta_\ell\} \subset \mathcal{B}_i$, or
between their reductions $\{\overline{\beta}_1, \dots,
\overline{\beta}_\ell\} \subset L_i$.

We now sketch how to find an explicit isomorphism of each $\mathcal{B}_i$
with a $\fieldk$-algebra of the form $\fieldk[x]/(h_i(x))$, in order to
obtain the isomorphism of~\eqref{equation2r.10}.  Finding such an
isomorphism is equivalent to finding a ``primitive element'' for the algebra
$\mathcal{B}_i$, which as we shall see is possible because $\fieldk$ is
perfect and because of the relation with the Dedekind domain $R$.  For
notational convenience, we shall drop the subscript $i$. 

\begin{proposition}
\label{proposition5r.2}
Given, as above, $\mathfrak{p} \subset \mathcal{B}$ with
$\dim_\fieldk \mathcal{B} = ef$, we can efficiently compute an element
$\beta \in \mathcal{B}$ whose minimal polynomial $h(x) \in \fieldk[x]$ has
degree $ef$.
\end{proposition}
\begin{proof}[Sketch of proof]
We first find a primitive element $\overline{\beta} \in L =
\mathcal{B}/\mathfrak{p}$, and its irreducible minimum polynomial
$g(x)\in \fieldk[x]$, where $\deg g(x) = f = [L:\fieldk]$.  This is
straightforward: for example, we can select random $\overline{\beta}$ (one
can show that the probability of selecting a primitive element is good),
and, for each candidate $\overline{\beta}$, find its minimal polynomial
$g(x)$ by looking for $\fieldk$-dependencies between  $\{1,
\overline{\beta}, \dots, \overline{\beta}{}^f\}$.  We repeat this process
until we find $\overline{\beta}$ for which $\deg g(x) = f$.  We now 
look for a lift $\beta \in \mathcal{B}$ of $\overline{\beta}$
whose minimal polynomial is $h(x) = (g(x))^e$.  This is trivial if $e=1$,
as any lift will do.  If $e \geq 2$, then we see that it suffices to
find a lift $\beta$ for which $g(\beta) \in \mathfrak{p} - \mathfrak{p}^2$
(since, in that case, $g(\beta) \in \mathcal{B}$ comes 
from an element of $R$ with valuation $1$ at the prime $\mathfrak{p}$).
Take an arbitrary lift $\beta_0$ of $\overline{\beta}$.  Since
$g(\overline{\beta}) = 0$, we know that $g(\beta_0) \in \mathfrak{p}$.  If
in fact $g(\beta_0) \notin \mathfrak{p}^2$, then we can take
$\beta = \beta_0$.  Otherwise, replace $\beta_0$ by $\beta_0 + \gamma$,
where we take any $\gamma \in \mathfrak{p} - \mathfrak{p}^2$.
This yields
\begin{equation}
\label{equation5r.7}
g(\beta_0 + \gamma) = g(\beta_0) + g'(\beta_0) \gamma + O(\gamma^2)
   \equiv g'(\beta_0) \gamma \pmod{\mathfrak{p}^2}.
\end{equation}
Since the extension $L/k$ is separable, we have $\g'(\overline{\beta}) \neq
0$, from which $g'(\beta_0)$ is a unit in $\mathcal{B}$, and we obtain what
we want.
\end{proof}






\providecommand{\bysame}{\leavevmode\hbox to3em{\hrulefill}\thinspace}
\providecommand{\MR}{\relax\ifhmode\unskip\space\fi MR }
\providecommand{\MRhref}[2]{%
  \href{http://www.ams.org/mathscinet-getitem?mr=#1}{#2}
}
\providecommand{\href}[2]{#2}


\end{document}